\documentclass{article}

\usepackage{PRIMEarxiv}
\usepackage{xcolor}

\usepackage{graphicx}
\usepackage{subcaption}
\usepackage{comment}
\usepackage{algcompatible}
\usepackage{mathtools}
\usepackage{booktabs}       
\DeclarePairedDelimiter{\floor}{\lfloor}{\rfloor}
\newcommand{\vect}[1]{\mathbf{#1}}

\usepackage{natbib}
 \bibpunct[, ]{(}{)}{,}{a}{}{,}%
 \def\bibfont{\small}%
\allowdisplaybreaks

\usepackage[utf8]{inputenc} 
\usepackage[T1]{fontenc}    
\usepackage{hyperref}       
\usepackage{url}            
\usepackage{booktabs}       
\usepackage{amsfonts}       
\usepackage{nicefrac}       
\usepackage{microtype}      
\usepackage{lipsum}
\usepackage{fancyhdr}       
\usepackage{graphicx}       
\graphicspath{{media/}}     

\title{Reinforcement Learning for Freight Booking Control Problems}

\author{
  Justin Dumouchelle \\
  University of Toronto \\
  \texttt{justin.dumouchelle@mail.utoronto.ca} \\
   \And
  Emma Frejinger \\
  Universit\'e de Montr\'eal  \\
  \texttt{emma.frejinger@umontreal.ca} \\
  \And
  Andrea Lodi \\
  Jacobs Technion-Cornell Institute,
  Cornell Tech and Technion - IIT \\
  \texttt{andrea.lodi@cornell.edu} \\
}

\pdfoutput=1

\begin{document}


\maketitle 

\begin{abstract}

Booking control problems are sequential decision-making problems that occur in the domain of revenue management.  More precisely, freight booking control focuses on the problem of deciding to accept or reject bookings: given a limited capacity, accept a booking request or reject it to reserve capacity for future bookings with potentially higher revenue. This problem can be formulated as a finite-horizon stochastic dynamic program, where accepting a set of requests results in a profit at the end of the booking period that depends on the cost of fulfilling the accepted bookings. For many freight applications, the cost of fulfilling requests is obtained by solving an operational decision-making problem, which often requires the solutions to mixed-integer linear programs.  Routinely solving such operational problems when deploying reinforcement learning algorithms that require large number of simulations may be too time consuming. The majority of booking control policies are obtained by solving problem-specific mathematical programming relaxations that are often non-trivial to generalize to new problems and, in some cases, provide quite crude approximations.  

In this work, we propose a two-phase approach: we first train a supervised learning model to predict the objective of the operational problem, and then we deploy the model within simulation-based reinforcement learning algorithms to compute control policies.  This approach is general: it can be used every time the objective function of the end-of-horizon operational problem can be predicted with reasonable accuracy, and it is particularly suitable to those cases where such problems are computationally hard. Furthermore, it allows one to leverage the recent advances in reinforcement learning as routinely solving the operational problem is replaced with a single prediction.  Our methodology is evaluated on two booking control problems in the literature, namely, distributional logistics and airline cargo management.

\end{abstract}


\section{Introduction }
\label{section:introduction}

The focus of this work is on quantity decisions arising in freight booking control problems. Such problems are of high practical relevance. This has been particularly manifest since the beginning of the pandemic, when we have experienced fast moving changes in freight transportation demand simultaneously as major disruptions in global supply chains occurred. This resulted in a highly uncertain environment with shifts in demand distributions over time. In this context, effective use of transport capacity is of the essence. Moreover, when faced with such a challenging environment, models and solution algorithms need to scale to large problems and should be flexible, i.e., it should be possible to adapt them regularly at a reasonable cost.

The booking control problem we consider has a relatively simple structure, and as we describe in the following, can straightforwardly be formulated as a Markov Decision Process (MDP).  We present a formulation similar to that delineated in \cite{talluri2006theory}. Let $\mathcal{N}$ denote a set of requests with cardinality $|\mathcal{N}|=n$.  The decision-making problem is defined over $T$ periods indexed by $t\in\{1,\ldots,T\}$. The probability that request $j\in\mathcal{N}$ is made at time $t$ is given by $p_j^t$ and the probability that no request occurs at period $t$ is given by $p_0^t$.  The time intervals are defined to be small enough such that at most one request is received in each period. 

To formulate this as an MDP, let $s_{j}^t$ denote the number of requests of type $j$ accepted before time $t$ and $\vect s_t \in \mathbb {R}^n, t=1,\ldots,T+1$, denote the state vector where the $j$-th index is given by $s_{j}^t$.  Furthermore, let $\mathcal{S}_t$ denote the set of all possible states at period $t$.  Since a given request can only be accepted or rejected, the actions are denoted by $a_j \in \{0,1\}$, where $a_j = 1$ if an offer for request $j$ is accepted.  The deterministic transition is given by $\vect s_{t+1} = \vect s_{t} + a_j\vect e_j$, with $\vect e_j \in \mathbb{R}^n$ having a $j$-th component equal to 1 and every other equal to 0. The revenue associated with accepting a request of type $j$ is given by $r_j$, and the reward function is given by $R: \mathcal{N} \times \{0, 1\} \rightarrow \mathbb{R}$, where $R(j, a_j) = r_j a_j$.  

An operational decision-making problem occurs at the end of the booking period and is related to the fulfillment of the accepted requests, henceforth referred to as the \emph{end-of-horizon problem} (EoHP). In the two applications we consider, they respectively correspond to  a capacitated vehicle routing problem (VRP) and a two-dimensional knapsack problem. The EoHP only depends on the state $\vect s_{T+1}$ and we denote its cost by $\Gamma: \mathbb{R}^n \rightarrow \mathbb{R}_{-}$.  With this, the value function is given by
\begin{equation}
    V_t(\vect s_t) = p_0^t V_{t+1}(\vect s_t) + \sum_{j\in\mathcal{N}} p_{j}^t \max_{a_j\in\{0,1\}} \{ r_j a_j + V_{t+1}(\vect s_t + a_j \vect e_j) \}, t = 1,\ldots, T,
    \label{eq:exact_mdp}
\end{equation}
\begin{equation}
    V_{T+1}(\vect s_{T+1}) = \Gamma(\vect s_{T+1}),
    \label{eq:exact_boundary}
\end{equation}
with $\vect s_1 = \vect 0$.

Due to the curse of dimensionality, solving this recursive formulation is intractable for problems of practical relevance \citep[e.g.,][]{adelman2007dynamic}. Furthermore, for reinforcement learning (RL) algorithms that typically require a large number of simulated system trajectories, the computational burden of learning a policy may be prohibitively high if hard EoHPs are solved routinely to compute (\ref{eq:exact_boundary}).  

The objective of our work is to design an approach that  (i) is broadly applicable to freight booking control problems with different EoHP structures, (ii) scales to large-size problems, and (iii) can learn from data. RL is appealing in this context, deep RL algorithms in particular as they have a high capacity, and are known to produce well-performing policies without an explicit model of system dynamics \citep[e.g.,][]{mnih2015human}. The main challenge associated with using RL for our problem class pertains to the EoHPs as they can be relatively costly to solve and good relaxations may not exist.  

In this paper, we precisely address this issue. In brief, we propose to embed fast predictions from a supervised learning algorithm tasked with handling (\ref{eq:exact_boundary}) within an RL approach. In particular, we focus our attention on the case where the EoHP is a mixed-integer linear program (MILP) and compare our approach against the state-of-the-art in two applications, namely, distributional logistics (DiL) introduced in \cite{ggm}, and airline cargo management (ACM), see e.g., \cite{barz2016air}.

Next, we provide an overview of the related work. We outline our contributions in Section~\ref{sec:contributions}, and close with a description of the remainder paper structure (Section~\ref{sec:paperStructure}).

\subsection{Related Work}

In this section, we discuss three areas of related work. First, we provide a brief background on the use of deterministic linear programming (DLP) to compute booking policies and we describe related approaches for solving the two booking control problems that we consider here. Second, we review existing work on  RL for booking control.  Third, we discuss machine learning (ML) for discrete optimization problems.

\subsubsection{Deterministic Linear Programming and Booking Control Problems.}
\label{sec:booking}

In revenue management (RM), DLP is by far the most widely adopted algorithmic approach for computing policies.  The primary reason for this is that DLP can be used to compute both booking-limit and bid-price policies through a linear program (LP) that directly models capacity constraints and the cost associated with fulfilling the requests.  The LP is formulated to maximize the quantity of each type of request that can be accepted without violating capacity constraints.  

Specifically, the DLP is given by
\begin{align}
\max & \quad \sum_{j\in\mathcal{N}} r_j \cdot x_j & \label{}\\
\text{subject to} & \quad A \vect x \le \vect b \label{},\\
            & \quad \vect 0 \leq \vect x \leq \vect \mu \label{},
\end{align}
where $\vect \mu \in\mathbb{R}^{n}$ is the vector of expected demand for each request type, $\vect A \in\mathbb{R}^{m\times n}$ and $\vect b \in\mathbb{R}^m$ denote the capacity constraints, and $\vect x \in\mathbb{R}^n$ is the variable that indicates the quantity of each type of request to be accepted.  Booking-limit control can be implemented directly through the primal solution by using the optimal solution ($\vect x^*$) as a threshold for the number of requests from class $j$ to accept for a single episode or trajectory.  Bid-price control can be imposed by requiring that the revenue from accepting a request is greater than a function of the dual variables and constraints.  This specific formulation does not account for the exact times and realizations of the requests. However, the LP can be re-solved throughout the booking period to provide updated primal and dual solutions.  Although DLP generally works well and is widely used in practice \cite[e.g.,][]{williamson1992airline,talluri2006theory}, it hinges on problem-specific assumptions and may be difficult to apply depending on the structure of the EoHP. Much of the early research in RM does not consider an EoHP and simply assumes a predetermined capacity cannot be violated. One of the first works considering an EoHP is \cite{bertsimas2003revenue}, where they propose dynamic programming and DLP-based approaches for simple EoHPs, namely overbooking and cancellations, in a general RM setting.  
As we discuss next, several studies have examined a broad diversity of specific problem settings and solutions methods.

ACM was introduced by \cite{kasilingam1997air} and has become a major focus in freight booking control over recent years.  Generally, the weight and volume of each request are uncertain until the departure.  At the end of the booking period, the set of accepted requests needs to be packed on an aircraft, and if there is insufficient capacity available, the unpacked items incur an off-loading penalty.  By only considering the weight and volume of the goods to be transported, the EoHP (technically a three-dimensional knapsack problem, see the discussion in Section \ref{sec:applications:cargo}) is formulated and solved as a vector-packing problem.  ACM can be formulated as a single-leg or a network problem, where the former focuses on control for a single flight and the latter is concerned with contexts in which cargo may need to be loaded on multiple flights to reach a final destination.  

\cite{amaruchkul2007single} propose a variety of heuristics based on DLP for single-leg ACM, where the off-loading penalty is simplified to a linear function of the excess weight and volume.
\cite{levin2012cargo} introduce a Lagrangian relaxation of the capacity constraints in the off-loading problem that are updated during the booking period to compute a policy.  
\cite{levina2011network} introduce an infinite horizon network cargo management problem with a solution based on approximate dynamic programming (ADP).
\cite{barz2016air} examine several approaches derived from DLP and ADP for single-leg and network problems while considering off-loading costs modeled as a vector packing problem (EoHP).

For DiL, the problem was only recently introduced by \cite{ggm}, and considers that a set of pick-up requests are received from a set of locations.  At the end of the booking period, the requests must be fulfilled, resulting in a capacitated VRP. This is a significantly more challenging EoHP as the exact solution to even a single VRP can be intractable. \cite{ggm} propose a MILP formulation, similar to a DLP formulation, but with the inclusion of integer variables to model the underlying routing problem.  A booking-limit policy is then given by the primal solutions of the MILP.  Unfortunately, solving the MILP exactly is intractable for even small instances, hence the approach in \cite{ggm} terminates after a small number of incumbent solutions are found in the solving process.

\subsubsection{Reinforcement Learning for Booking Control.} \label{sec:reviewRL}

In recent years, RL has been successful in a variety of challenging tasks \citep{silver2016mastering, mnih2015human}. However, its application is currently foremost restricted to simulation environments (such as games), or real-world settings where system failure is not costly and of limited risk \citep{donti2021enforcing}. As we discuss below, its application is still limited within RM, in particular in the context of freight applications. There is an extensive and fast-evolving literature on RL so a comprehensive overview is out of the scope of this paper. Instead, we refer the reader to excellent textbooks \citep{bertsekas2019reinforcement,sutton2018reinforcement} and surveys \citep[e.g.,][]{Gosavi09, BusonEtAl18}.

Within RM, RL has been applied to solve pricing problems under various market assumptions \citep[e.g.,][]{KephTesa00, RanaEtAl14, ShuklaEtAl20,  KastSchl22}, and to airline seat inventory control problems \citep[e.g.,][]{gosavii2002reinforcement,lawhead2019bounded,bondoux2020reinforcement,shihab2019autonomous}. Both contexts have distinguishing features compared to freight booking control problems. Importantly, one passenger occupies exactly one seat and the capacity is defined as a given number of seats. On the contrary, in the airline cargo and DiL applications discussed above, deploying RL is more involved because of the potentially hard EoHP. 

We are only aware of one study that uses RL for freight booking control in a realistic setting. \cite{SeoEtAl21} consider a sea cargo booking control problem. Instead of profit maximization, they maximize revenue, hence ignoring both the cost associated with transporting accepted demand and constraints related to loading the cargo. In other words, they ignore the EoHP and RL can be straightforwardly applied.

\subsubsection{Machine Learning for Discrete Optimization.}
\label{sec:ML4DO}

The relatively recent success of ML in a wide range of applications has led to a fast-growing interest within the optimization community. In particular, the area at the intersection of ML and discrete optimization has been expanding from an impressive surge of contributions \citep[see, e.g.,][for a methodological analysis of the literature]{bengio2020machine}. 

ML approaches for discrete optimization problems have been developed to either predict the entire solution of a problem, i.e., replacing the discrete optimization techniques by ML ones, or to accelerate discrete optimization methods \citep[referred to as ML-augmented combinatorial optimization in the survey][]{KotaryEtAl21}. We refer to \cite{bengio2020machine} for references on these topics. Instead, we concentrate on using ML to predict the objective value of a (hard) discrete optimization problem, a task that has received less attention with a few notable exceptions:  \cite{fischetti2017using} use supervised learning to predict the optimal objective value of MILP problems for a wind energy application. \cite{prates2019learning} propose a combination of supervised learning and binary search to predict optimal solution values to the traveling salesperson problem (TSP). \cite{LarsEtAl22} focus on predicting aggregate solutions (objective function value is a special case) to stochastic discrete optimization problems in the context of two-stage stochastic programs.

Closest to our work are heuristic methods that integrate predictions from supervised learning to gain speed. \cite{LarsFrejGendLodi22} propose a matheuristic version of a Benders decomposition method for solving two-stage stochastic programs. \cite{dumouchelle2022neur2sp} propose easy-to-solve surrogate models for two-stage stochastic programs based on neural networks trained to predict the second-stage objective.  For multi-stage stochastic programming, \cite{dai2022neural} propose a supervised approach for learning piecewise linear convex value function approximations.  Whereas the combination of supervised learning and RL is not new \citep[e.g.,][]{ye2003fuzzy, henderson2005hybrid, finn2017generalizing}, we are unaware of any works that introduce predictions of optimal discrete optimization problem values within an RL algorithm.

\subsection{Contributions} \label{sec:contributions}
In this work, we propose a two-phase approach, wherein we first train a supervised learning model to predict the value of the EoHP solution, and then we integrate predictions within simulation-based RL algorithms to compute control policies. This approach is general as it can be applied to any EoHP, regardless of its structure, as long as the objective value can be computed or approximated. In addition, as the prediction significantly reduces the time required for simulation, this can be scaled to complex EoHPs and allow for the use of RL algorithms requiring extensive simulation.  

We make the following contributions:
\begin{enumerate}
    \item We propose a methodology that (i) formulates an approximate MDP by replacing the (discrete) optimization model of the EoHP with a prediction function yielded by offline supervised learning, and (ii) uses RL techniques to obtain approximate control policies, both of which can be implemented with existing off-the-shelf libraries.
    \item We apply the proposed methodology to DiL and ACM problems from the literature and we show computationally that we achieve high-quality control policies that can be deployed efficiently. 
    The EoHP is harder for the DiL application. In this case, the percent improvement in mean profit compared to the algorithms presented in \cite{ggm} is at least 28.6\% 
    with an average online computing time of less than 4.5 seconds.
    Notably, we achieve a performance close to the strong state-of-the-art baseline for the well-studied ACM problems. In this case, the percent improvement in mean profit compared to the baseline method ranges between -2.2\% (baseline is better) and 1.1\% with an average online computing time of less than 1.5 seconds.
    In all cases, the total offline computing time is less than six hours, orders of magnitude smaller than it would be without leveraging predictions from a supervised learning model.
\end{enumerate} 

Our work is of high practical importance as it can deal with EoHPs of varying structure using existing off-the-shelf libraries for supervised and reinforcement learning. In contrast, existing work rely on simplifications of the EoHP or tailor-made methodologies leveraging specific problem structure. Our methodology is relatively easy to use in practice for EoHP that have a different structure compared to those treated in the literature.

\subsection{Structure of the Paper} \label{sec:paperStructure}

The remainder of the paper is organized as follows. In Section~\ref{section:methodology}, we introduce our methodology.  
In Section~\ref{section:applications}, we discuss two freight booking control problems, their state-of-the-art baselines, and how our supervised prediction task is formulated.
In Section~\ref{section:results}, 
we compare our methodology to the problem-specific baselines.
Finally, Section~\ref{section:conclusion} draws conclusions and highlights avenues for future work.
\section{Methodology}
\label{section:methodology}
We divide the description of our methodology into three parts: In Section~
\ref{sec:APPROX}, we propose an approximation of the MDP \eqref{eq:exact_mdp}--\eqref{eq:exact_boundary}, and, in Section \ref{sec:EOH}, we specify the ML ingredients of that approximation. Finally, in Section \ref{sec:RL:motivation}, we motivate the use of RL to solve the approximate MDP of Section \ref{sec:APPROX}.

\subsection{Approximate MDP Formulation}
\label{sec:APPROX}

Solving the MDP \eqref{eq:exact_mdp}--\eqref{eq:exact_boundary} can be intractable even for small instances. A possible approach to deal with this is to approximate \eqref{eq:exact_boundary}. This is especially important in the context of booking control problems, where it may be costly to compute. Specifically, the approximation of $\Gamma$ is denoted as $\phi : \mathbb{R}^m \rightarrow \mathbb{R}_-$ and given by
\begin{equation}
    \Gamma(\vect s_{T+1}) \approx \phi(g(\vect s_{T+1})),
    \label{eq:approx_operational}
\end{equation}
where $g:\mathbb{R}^n \rightarrow \mathbb{R}^m$ is a mapping from the state at the end of the horizon to an $m$-dimensional representation. The approximate MDP is then given by
\begin{equation}
    \tilde V_t(\vect s_t) = p_0^t \tilde V_{t+1}(\vect s_t) + \sum_{j\in\mathcal{N}} p_{j}^t \max_{a_j\in\{0,1\}} \{ r_j a_j + \tilde V_{t+1}(\vect s_t + a_j \vect e_j) \}, t = 1, \ldots, T,
    \label{eq:approx_mdp}
\end{equation}
\begin{equation}
    \tilde V_{T+1}(\vect s_{T+1}) = \phi(g(\vect s_{T+1})).
    \label{eq:approx_boundary}
\end{equation}

Note that in the system \eqref{eq:approx_mdp}--\eqref{eq:approx_boundary}, the only change compared to \eqref{eq:exact_mdp}--\eqref{eq:exact_boundary} is in the value function at the end of the time horizon. However, the value will be recursively propagated to the entire value function.

As we further motivate in Section~\ref{sec:RL:motivation}, we focus on solving the approximate MDP in \eqref{eq:approx_mdp}-\eqref{eq:approx_boundary} with RL. Such algorithms are based on fast evaluations of system trajectories. In turn, this requires evaluating \eqref{eq:approx_boundary} and a related fast approximation scheme is the topic of the following section.

\subsection{Approximating the End-of-Horizon Cost}
\label{sec:EOH}
\label{section:methodology:supervised}

The approximation of $\Gamma(\cdot)$ in~\eqref{eq:approx_boundary} can be computed in various ways, for instance by using a problem-specific heuristic, solving a MILP (or its LP relaxation) to a given time limit or optimality gap or, as proposed in this paper, it can be predicted by ML. 
In this work, we focus our attention on prediction through supervised learning, i.e., we develop a mapping from the end-of-horizon state to the associated operational cost. Specifically, a supervised learning model is trained offline to separate the problem of accurately predicting $\Gamma(\cdot)$ from that of solving \eqref{eq:approx_mdp}--\eqref{eq:approx_boundary}.  

The choice of using supervised learning is motivated by the three following considerations.  First, predicting the EoHP cost is problem agnostic as long as a set of features and labels can be constructed. This allows our methodology to be easily adapted to new applications, which is not the case for problem-specific heuristics, where minor changes in the problem formulation often require non-trivial -- if at all possible -- adaptations to be made to the heuristic.  
Second,  the use of supervised learning is motivated by efficiency considerations.  As we only require a single prediction for each simulated trajectory, this can be made in very short computing time (milliseconds), often orders of magnitude faster than heuristics for mathematical programming problems, especially when considering challenging EoHP. Third, the sample size required to learn the EoHP cost separately from the booking policy is likely much smaller than that required if they were learned jointly. We analyze this numerically in Section~\ref{section:results}.

The aim of the supervised learning algorithm is to accurately predict the value of the EoHP for any feasible (and relevant) state $\mathbf{s}_{T+1}$. That is, we aim to achieve accurate generalization for states unseen in the training data. For this purpose, we use simulation to collect $N$ labeled data points $\mathcal{D} = \{(\vect s_{T+1}^1, \Gamma(\vect s_{T+1}^1)), \ldots, (\vect s_{T+1}^N, \Gamma(\vect s_{T+1}^N))\}$ and an associated feature mapping $g(\cdot)$. Next, we describe the related steps in detail.

\paragraph{Sampling End-of-Horizon States $\mathbf{s}_{T+1}$.} We use a stationary random policy that accepts a request with a given probability $p$ to simulate trajectories.  At the end of the time horizon, a state $\vect s_{T+1}$ is obtained. 
This process is repeated for different values of $p$. The objective is to generate a distribution of feasible final states in the data (optimal and sub-optimal ones) that will lead to high-quality generalization of the ML algorithm.

\paragraph{Feature Mapping.} 

The feature mapping is problem specific. It is represented by the function $g(\cdot)$ in~\eqref{eq:approx_boundary}. In general, for booking-control problems, we will have accepted a \emph{set} of requests. Since a set is a permutation-invariant structure, it is desirable to use a permutation-invariant learning model similar to those described in \cite{zaheer2017deep} to capture the features of the set of accepted requests. In short, this type of learning model achieves permutation invariance with the use of shared parameters for each request and a permutation-invariant aggregation operation that combines the learned representations of all accepted requests.  In contrast, other learning models require aggregating request features to achieve permutation invariance, resulting in a loss of precise request information, which may be crucial for accurate prediction.

\begin{figure}[h!]
    \centering
    \includegraphics[width=1.0\textwidth]{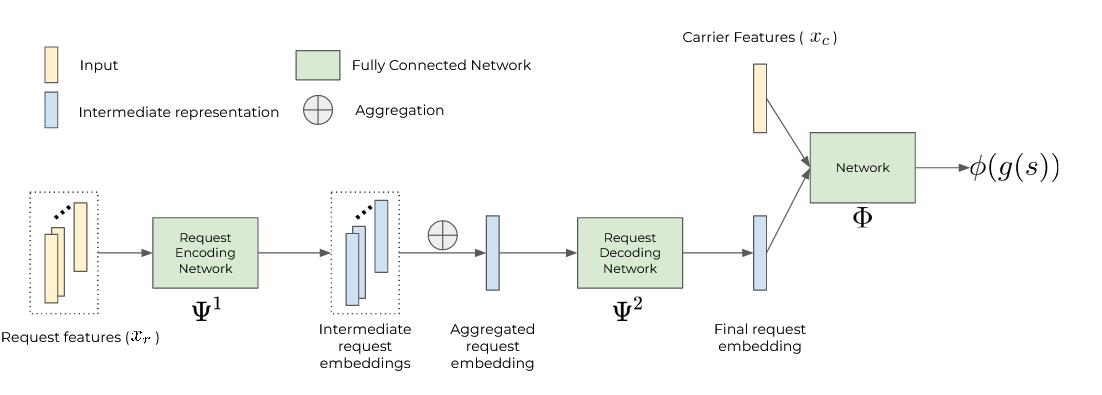}
    \caption[Pictorial representation of permutation-invariant model for booking control. ]%
    {{\small Pictorial representation of permutation-invariant model for booking control.  Adapted from Figure 2 in \cite{dumouchelle2022neur2sp}}.}    
    \label{fig:nn_arch}
\end{figure}

For the permutation-invariant modeling (see~Figure~\ref{fig:nn_arch} for an illustration), the vectors of features, say $x_{r}$, characterizing the individual requests are independently passed through an encoding feed-forward network, $\Psi_1$. The outputs of this network are then aggregated with a permutation-invariant transformation (a sum in our case). The resulting permutation-invariant output is passed through a decoding feed-forward network, $\Psi_2$, thus yielding the permutation-invariant embedding vector of all requests, as advocated in \cite{zaheer2017deep}. Next, the concatenation of this embedding vector with the vector characterizing the features of the carrying vessels and global information, say $x_{c}$, is passed through a third feed-forward network, $\Phi$, yielding a prediction for the EoHP cost.  

As an alternative, one could replace the embedding with a vector that directly aggregates the individual statistics of the accepted requests. However, as evidenced in Section~\ref{sec:SL}, models that use the aggregated features suffer from notably high prediction errors.

\paragraph{Computing Labels.} For each sampled end-of-horizon state $\mathbf{s}_{T+1}^n,~n=1,\ldots,N$, we compute labels $\Gamma(\vect s_{T+1}^n)$ by solving the EoHP. This can be done with an exact method, or a heuristic depending on the required precision and computing time budget.  

An important aspect of this data generation procedure is that solving the EoHP can be parallelized up to the number of available threads.  This stands in contrast with online approaches for estimating the EoHP cost (such as a problem-specific heuristic or LP-relaxation) as they involve significantly more costly simulation time when training RL algorithms.  In addition, since most RL algorithms rely on sequential simulation, parallelization is usually more difficult.  

\subsection{Policies for the Approximate MDP}
\label{sec:RL:motivation}

Once supplied with an offline supervised learning model trained to approximate the EoHP cost, we can leverage a broad range of RL techniques to compute control policies by solving~\eqref{eq:approx_mdp}-\eqref{eq:approx_boundary}.
We select the Q-learning algorithm \cite{watkins1992q}. Its performance has been well documented through extensive applications on a broad range of tasks and it is known to feature typically short computing times.
More specifically, we focus on Deep Q-learning (DQN).

We implement two variants of DQN.  The first is DQN-L that simply considers a linear state representation containing the number of accepted requests for each class along with additional information about the utilized capacity and period.  DQN-L employs a fully connected network to approximate the value function.  Similar to the supervised learning case, we have a set of accepted requests during simulation, so once again permutation-invariant approaches are a natural choice for the DQN value-function estimator.  This motivates the choice of the second model, DQN-S, which is similar in architecture to Figure~\ref{fig:nn_arch} as it accounts for the features of each of the accepted requests from previous periods.

\section{Applications}

\label{section:applications}

In this section, we introduce two specific booking control problems, describe their respective state-of-the-art baselines, and examine considerations pertaining to the prediction of the EoHP costs.

\subsection{A Distributional Logistics Problem}
\label{sec:applications:dl}

We start by describing the formulation of the EoHP proposed by \cite{ggm} for a DiL problem. We then outline their solution approach and follow with a description of the associated supervised prediction task (Section~\ref{sec:LogisticsPredTask}). We close with a discussion of RL training in Section~\ref{sec:LogisticRLTraining}.

\subsubsection{EoHP Formulation}

In this problem, booking requests correspond to pickup activities and the EoHP to a capacitated VRP. Each pickup request has an associated location and revenue. The cost incurred at the end of the booking period is the cost of the VRP solution and hence depends on all the accepted requests.

The problem can be formulated as \eqref{eq:exact_mdp}--\eqref{eq:exact_boundary} with the boundary function described below.
The end-of-horizon VRP consists of a fixed set of $K_0$ vehicles, each with capacity $Q\in\mathbb{R}_+$.
The depot location is indexed by location $0$. Each request is associated with a location, so we use $\mathcal{N}$ to denote the set of request locations. The set of all nodes $\mathcal{V}$ is given by $\mathcal{N} \cup \{0\}$, i.e., the union of the request locations and the depot.  The set of directed arcs is given by $\mathcal{V} \times \mathcal{V}$ and denoted as $\mathcal{A}$.  
The cost of an arc is given by $c_{ij}\in\mathbb{R}_+$, for $(i,j)\in\mathcal{A}$.
The optimal objective value is denoted by $z^*(\vect s, K)$, where $K$ denotes the number of vehicles used.  
If more than $K_0$ vehicles are required, then additional outsourcing vehicles are allowed to be used at an additional fixed cost of $C\in\mathbb{R}_+$ each.
\cite{ggm} do not consider the inclusion of outsourcing vehicles and instead assign an infinite operational cost if the requests cannot be fulfilled.  
To accomplish the same objective, $C$ can be chosen sufficiently large enough such that the cost of adding an additional vehicle is larger than any potential revenue from the requests it can fulfill. The operational cost is given by

\begin{equation}
\label{eq:vrp_cost}
\Gamma(\vect s) = \max_{K \ge K_0, K\in\mathbb{Z}_+} \{- z^*(\vect s, K) - C(K-K_0)\},
\end{equation}
where $z^*(\vect s, K)$ denotes the cost of the routing problem given $K$ vehicles, and $C(K-K_0)$ denotes the cumulative cost of the $K-K_0$ outsourcing vehicles used.  The routing problem can be represented as a MILP and, following \cite{ggm}, is given by

\begin{align}
   z^*(\vect s, K)  = &\min\quad \sum_{k=1}^{K} \sum_{(i,j)\in\mathcal{A}} c_{ij}\alpha_{ij}^k      \label{eq:vrp_obj} \\[2.5ex]
  \text{s.t.\quad}  & \sum_{(i,j)\in\mathcal{A}} \alpha_{ij}^k = \beta_i^k       &  & \forall i\in\mathcal{V}, \forall k\in\{1,\ldots,K\}     \label{eq:vrp_arc_in},\\[1.5ex]
                    & \sum_{(i,j)\in\mathcal{A}} \alpha_{ji}^k = \beta_i^k       &   & \forall i\in\mathcal{V}, \forall k\in\{1,\ldots,K\} \label{eq:vrp_arc_out}  ,\\[1.5ex]
                    & \sum_{k=1}^{K} \beta_{j}^k \le 1                           &   & \forall j\in\mathcal{N}  \label{eq:vrp_visited_once},\\[1.5ex]
                    & \sum_{k=1}^{K} \beta_{0}^k \le  K ,                        &   & \label{eq:vrp_num_routes} \\[1.5ex]
                    & \sum_{i\in \mathcal{S}}\sum_{j\in\mathcal{V}\backslash\mathcal{S}}  \alpha_{ij}^k \ge \beta_h^k       
                                                                                 &   & \forall \mathcal{S}\subset \mathcal{V}: 0\in\mathcal{S}, \forall h\in\mathcal{V}\backslash\mathbb{S}, \forall k\in\{1,\ldots,K\} \label{eq:vrp_sec} ,\\[1.5ex]
                    & q_j^k \le Q \beta_j^k                           &   & \forall j\in\mathcal{N}, k\in\{1,\ldots,K\} \label{eq:vrp_single_capacity} ,\\[1.5ex]
                    & \sum_{j\in\mathcal{N}} q_j^k \le Q            &   & \forall k\in\{1,\ldots,K\} \label{eq:vrp_capacity} ,\\[1.5ex]
                    & \sum_{k=1}^{K} q_j^k = s_j                    &   & \forall j\in\mathcal{N} \label{eq:vrp_fulfilled} ,\\[1.5ex]
                    & \alpha_{ij}^k\in\{0,1\}                       &   & \forall (i,j)\in\mathcal{A}, \forall k\in\{1,\ldots,K\} \label{eq:vrp_alpha} ,\\[1.5ex]
                    & \beta_i^k \in\{0,1\}                          &   & \forall i\in\mathcal{V}, \forall k\in\{1,\ldots,K\} \label{eq:vrp_beta} ,\\[1.5ex]
                    &  q_j^k \ge 0                  &   & \forall j\in\mathcal{N}, \forall k\in\{1,\ldots,K\}. \label{eq:vrp_q}
\end{align}

The decision variables are $\alpha_{ij}^k$, $\beta_i^k$, and $q_j^k$:  The binary variable $\alpha_{ij}^k$ equals $1$ if vehicle $k$ uses arc $(i,j)$.  The binary variable $\beta_i^k$ equals $1$ if vehicle $k$ visits location $i$.  The quantity of accepted requests from vehicle $k$ at location $j$ is given by the non-negative variable $q_j^k$.  Objective function~\eqref{eq:vrp_obj} minimizes the routing cost. Constraints~\eqref{eq:vrp_arc_in} and~\eqref{eq:vrp_arc_out} ensure that each vehicle goes in and out of the visited locations. Constraints~\eqref{eq:vrp_visited_once} assert that each location must be visited at most once, and \eqref{eq:vrp_num_routes} limits the number of vehicles to at most $K$. Constraints \eqref{eq:vrp_sec} ensure that every route traveled by a vehicle is connected, and each constraint in \eqref{eq:vrp_single_capacity} restricts vehicle $k$ from collecting anything from location $j$ if it does not visit the location (i.e., if $\beta_j^k = 0)$. Finally, constraints~\eqref{eq:vrp_capacity} guarantee that the capacity $Q$ of each vehicle is not exceeded, and \eqref{eq:vrp_fulfilled} ensure that the accepted requests at each location are fulfilled. It is well known that VRP formulations with arc variables duplicated for each vehicle, i.e., $\alpha_{ij}^k$ in \eqref{eq:vrp_obj}--\eqref{eq:vrp_q} above, are not the most efficient ones \citep{toth2002vehicle}. However, this formulation is suitable for dealing with the booking control version of the problem \citep[see][for details]{ggm}. 

\subsubsection{State-of-the-art Solution Method}
\label{sec:VRP:SOTA}

The DLP-based approach of \cite{ggm} constitutes the state-of-the-art for the DiL booking control problem.
They formulate a profit-maximizing MILP which is given by
\begin{align}
   \max\quad  & \sum_{j\in\mathcal{N}} r_j x_j - \sum_{k=1}^{K} \sum_{(i,j)\in\mathcal{A}} c_{ij}\alpha_{ij}^k      \label{eq:pmvrp_obj} \\[2.5ex]
  \text{s.t.\quad}  & \eqref{eq:vrp_arc_in} - \eqref{eq:vrp_capacity}, \eqref{eq:vrp_alpha} - \eqref{eq:vrp_q}, \\
    & \sum_{k=1}^{K} q_j^k = s_j + x_j &   \forall j\in \mathcal{N} &  \label{eq:pmvrp_demand}, \\
    & 0 \le x_j \le \mu_j  &   \forall j\in\mathcal{N} &   ,\label{eq:pmvrp_y} 
\end{align}

where $x_j$ is a variable indicating the number of requests to accept at location $j$, $\mu_j$ is the expected number of requests for location $j$, and $s_j$ represents the number of previously accepted requests for location $j$.
The Booking Limit Policy (BLP) simply solves \eqref{eq:pmvrp_obj}-\eqref{eq:pmvrp_y} before any requests have been accepted, i.e., $s_j=0, \forall j\in\mathcal{N}$.  A policy is then given by defining thresholds for each location given by the primal solution $x_j^*$.  The Booking Limit Policy with Reoptimization (BLPR) is similar to BLP. However, the threshold policy is updated halfway through the booking period.  This is achieved by updating $s_j$ based on the number of accepted requests in the current trajectory/episode, then the threshold is defined from the new MILP solution $\vect x^*$.  These approaches are then applied by accepting requests until the threshold for a given class is met. Because the DLP-type formulation is not, in fact, a LP but a MILP, the dual solution is not available, so it cannot be used to compute bid-price policies.  

Considering the complexity of the EoHP, solving the MILP is difficult and \cite{ggm} limit the process to two incumbent solutions for the profit maximization problem and four incumbent solutions for the end-of-horizon MILP. Essentially, the MILP run is aborted after two (or four) primal solutions have been found, independently of their quality. This formulation additionally suffers from the fact that low probability-high reward requests may be rejected if their expected demand is quite small, which is perhaps the largest limitation of the BLP/BLPR in this context.

\subsubsection{Prediction Task} \label{sec:LogisticsPredTask}

The objective is to produce an accurate approximation of~\eqref{eq:vrp_cost} that is fast to compute.  In general, one can predict the sum of the MILP objective function value and the outsourcing cost, i.e.~\eqref{eq:vrp_cost}. 
However, this is likely to be detrimental as the outsourcing cost may result in a large change in the label whereas only a minor change occurs in the features.
For this reason,  ML is used to predict $z^*(\vect s, K)$ while the outsourcing cost $C(K-K_0)$ is computed separately by solving the associated bin-packing problem of determining the minimum number of vehicles needed to do all pickup activities. The bin-packing is solved to optimality by the solver MTP \citep{martello1990knapsack}.  

For this DiL application, even solving the EoHP to optimality to collect labeled data is prohibitively expensive. For this reason, we leverage a VRP-specific heuristic that allows us to compute high-quality solutions within reasonable computing times. 
Specifically, we use fast iterated local search localized optimization (FILO),  \citep{accorsi2021fast} to compute a label \eqref{eq:vrp_obj}.
One minor consideration we had to make when using FILO was to ensure the number of vehicles was minimized, which is not done explicitly in FILO as it is an unlimited fleet solver.  To bias FILO, we simply add a prohibitively expensive cost for leaving the depot, which in turn minimizes the number of vehicles leaving the depot.  The additional cost is then subtracted after a solution is computed. 

We use a permutation-invariant model and compare it to a fixed-size input model. For the permutation-invariant model, the request features, $x_{r}$, are given by the location index, x and y coordinates, and demand.  The carrier features, $x_{c}$, are given by the number of vehicles, the number of locations with demand, the x and y coordinates of the depot, and the capacity. 
For the directly aggregated input structure, the
capacity, depot location, total number of accepted requests per location and aggregate statistics of the locations are computed, namely, the distance between locations and the depot, and relative distances between locations are used.  For each of these, we compute the min, max, mean, median, standard deviation, and 1st/3rd quartiles.

\subsubsection{RL Training} \label{sec:LogisticRLTraining}

For RL, we have two distinct state space representations.  For DQN-L, we consider a linear state space with the following features: the incoming request index, the number of accepted requests at each location, and the period.  For DQN-S, we consider a permutation-invariant model similar to that in Section~\ref{section:methodology:supervised} with request features, $x_{r}$, given by the location index and the x and y coordinates.  The carrier features, $x_{c}$, are given by the total capacity, utilized capacity, number of vehicles, and the period. The time period can be indicated through either an integer or a one-hot vector. Here, one-hot refers to the vector the size of the maximal number of periods with a 1 in the index of the current period and 0's in all other indices.  One of these hyperparameter settings is selected at the time of model selection.

\subsection{An Airline Cargo Problem}
\label{sec:applications:cargo}

We introduce the ACM problem described in \cite{barz2016air} \citep[originally proposed by][]{amaruchkul2007single}, and discuss the details required for building a predictor for the EoHP.  
As briefly mentioned, an ACM problem can be formulated as either a single-leg or a network problem.
A single-leg problem is limited to cargo traveling between an origin-destination pair, while a network problem allows for cargo to travel on multiple flights within a network.  In this work, we focus on the single-leg problem as it is more rigorously studied and the generalization to the network problem is straightforward with the methodology presented here.  
The MILP formulation of the EoHP available in \cite{barz2016air} can be specialized to the single-leg context as follows.

\subsubsection{EoHP Formulation}
In the single-leg cargo-management application, booking requests correspond to shipments that are being transported on commercial airline flights. Typically, commercial airlines have a non-negligible quantity of space remaining after passenger luggages are stored.  With the excess capacity, commercial shipments are often transported between origin and destination, resulting in a significant source of additional revenue.  
Each request is considered to be of a particular class, such as electronics, fresh food, or raw materials, which share a joint weight and volume distribution.

The EoHP incurred at the end of the booking period is given by a vector packing problem 
\citep[see, e.g.,][]{caprara2001lower}.
The EoHP cost is given by the set of accepted requests that cannot be loaded due to the capacity of the aircraft and is referred to as the off-loading cost. More precisely,
each request has an associated weight and volume and the flight has both a weight and volume capacity. 
Whereas all of these quantities are uncertain until the time of departure, their expected values are assumed known.

We now define the MILP formulation for the EoHP.  The realized weight and volume capacities of the flight are given by $Q_w\in\mathbb{R}_+$ and $Q_v\in\mathbb{R}_+$, respectively. The set of shipping classes or request types is given by $\mathcal{N}$. The realized weight and volume associated with the $i$-th request of the $j$-th shipping class are denoted by $w_{ji}$ and $v_{ji}$, respectively.  For airline cargo, it is often convenient to express the profit and off-loading costs as a function of the weight and volume.  
This is done through the concept of \emph{chargeable weight}, which is defined by a weighted maximum of the weight and volume for a given request.  The chargeable weight allows requests to be charged in terms of the predominant capacity dimension. 
Specifically, the chargeable weight is given by $\max\{w_{ji}, v_{ij} / \eta \}$, where $\eta\in\mathbb{R}_+$ is a constant representing the weight-volume ratio of a shipment.  Using the chargeable weight, the off-loading cost is then given by $c_{j}\cdot \max\{w_{ji}, v_{ij} / \eta \}$.  Finally, the operational cost is given by
\begin{equation}
\label{eq:vp_cost}
\Gamma(\vect s) = - z^*(\vect s),
\end{equation}
where $z^*(\vect s)$ denotes the off-loading cost associated with the vector packing problem
\begin{align}
   z^*(\vect s)  = &\min\quad \sum_{j\in\mathcal{N}}
   \sum_{i=1}^{s_j} c_j\max\{ w_{ji}, v_{ji} / \eta \} (1 - y_{ji})    \label{eq:vp_obj} \\[2.5ex]
  \text{s.t.\quad}  &
  \sum_{j\in\mathcal{N}}\sum_{i=1}^{s_j} y_{ji} w_{ji} \le Q_w   , &  & \label{eq:vp_wei}\\[1.5ex]
  & \sum_{j\in\mathcal{N}}\sum_{i=1}^{s_j} y_{ji} v_{ji} \le Q_v  ,&  & \label{eq:vp_vol}\\[1.5ex]
  &  y_{ji} \in\{0,1\}                & \quad  \forall j\in\mathcal{N}, \forall i\in\{1,\ldots,s_j\}.  & \label{eq:vp_var}
\end{align}
The binary variables $y_{ji}$ are equal to $1$ if the $i$-th accepted request of the $j$-th shipping class is loaded, and $0$ otherwise.  
Objective function~\eqref{eq:vp_obj} minimizes the sum of the off-loading costs of the shipments multiplied by a constant $c_j$ associated with each shipping class $j$. Constraints~\eqref{eq:vp_wei} and~\eqref{eq:vp_vol} ensure that the weight and volume of the loaded items respect the capacity of the flight.

Although modeling the loading requirements as the two independent knapsack constraints on weight and volume \eqref{eq:vp_wei}--\eqref{eq:vp_vol} is an established practice in the literature, it is worth mentioning that this is clearly a crude approximation. In fact, the real loading problem would require to consider three-dimensional objects (instead of approximations with their volume) together with their weight  \citep[see, e.g., the extensive discussion in][]{nickel}. The corresponding optimization problem would be a generalization of the already very complex three-dimensional knapsack. 

\subsubsection{State-of-the-art Solution Method}
\label{sec:CM:SOTA}

For ACM, the set of available algorithms is significantly larger than for the DiL problem.  \cite{barz2016air} provide an in-depth numerical analysis of algorithms, but the top-performing algorithms are all based on the DLP.  As the weight and volume are unknown until the end of the booking period, the DLP is formulated using their expectations.  Specifically, the expected weight and volume capacities are denoted by $E[Q_w]$ and $E[Q_v]$, respectively.  The expected weight and volume of an item of class $j$ are denoted by $E[w_j]$ and $E[v_j]$, respectively.  With these assumptions and definitions, the DLP is then given by 

\begin{align}
   \max\quad & \sum_{j\in\mathcal{N}} r_j x_j - c_j z_j
   \label{eq:pmvp_obj} \\[2.5ex]
  \text{s.t.\quad}  &
    \sum_{j\in\mathcal{N}} y_{j} \le E[Q_w],        &  & \label{eq:pmvp_wei}\\[1.5ex]
    & \sum_{j\in\mathcal{N}} y_{j} \le E[Q_v],      &  & \label{eq:pmvp_vol}\\[1.5ex]
    & z_j + y_j \ge x_j E[w_j]                                  & \forall j\in\mathcal{N} & \label{eq:pmvp_wei_res},\\[1.5ex]
    & z_j + \frac{y_j}{\eta} \ge x_j \frac{E[v_j]}{\eta}          & \forall j\in\mathcal{N} & \label{eq:pmvp_vol_res},\\[1.5ex]
    &  y_{j}, z_j \ge 0                                         & \quad  \forall j\in\mathcal{N} & \label{eq:pmvp_var_xz},\\[1.5ex]
    &  0 \le x_{j} \le \mu_j                                    & \quad  \forall j\in\mathcal{N} & \label{eq:pmvp_var_y},
\end{align}
where $x_j$ models the accepted requests for a class $j$. Variables $y_j$ and $z_j$ model the expected weight and volume of a class with $y_j$ representing the quantity that can be loaded and $z_j$ representing the quantity that cannot.  The objective \eqref{eq:pmvp_obj} models the profit from accepting requests, minus the offloading costs.  Constraints \eqref{eq:pmvp_wei}-\eqref{eq:pmvp_vol} ensure that the capacity is not exceeded.  Constraints \eqref{eq:pmvp_wei_res}-\eqref{eq:pmvp_vol_res} model the amount that is loaded/offloaded.  Finally, constraints \eqref{eq:pmvp_var_xz}-\eqref{eq:pmvp_var_y} define the variable restrictions.

It is important to note that this formulation is an LP relaxation to the original vector packing problem. However, since the LP relaxation of the vector packing is relatively tight \citep{caprara2001lower}, this will likely not have a large impact on the objective value.  Furthermore, since we have an LP, duality can be used to compute bid-price policies.  In \cite{barz2016air}, the DLP algorithm is defined by using the dual variables of constraints \eqref{eq:pmvp_wei}-\eqref{eq:pmvp_vol} to model opportunity cost from accepting a given request.  Several algorithms are then built on the idea of using the DLP solution, with the best-performing algorithm being dynamic programming decomposition (DPD).  This algorithm decomposes the weight and volume into a single quantity using the dual variables, then formulates a tractable one-dimensional MDP that can be solved via backtracking \cite[see][for details]{barz2016air}.

The practicality of this approach is highly dependent on the simplification of the EoHP discussed in the previous section. Solving the real generalized three-dimensional knapsack would have two major consequences: a computationally much more demanding problem to be solved repeatedly, and the loss of duality to derive bid-price policies because  
the LP relaxation of packing problems whose constraints are dependent on 3-dimensional objects is extremely poor.

\subsubsection{Prediction Task}
The objective is to produce an accurate approximation of~\eqref{eq:vp_cost} that is fast to compute. For this purpose, ML is used to predict $z^*(\vect s)$.  In contrast to the previous application where the EoHP required a heuristic to label the data, for airline cargo the EoHP can be solved exactly within a reasonable computing time, so we use the optimal objective in this case.  

Similar to the first application, we use a permutation-invariant model as well as a model with fixed-size input.  
For the permutation-invariant model, the request features, $x_{r}$, are given by the weight, volume, and cost of each request.  The carrier features, $x_{c}$, are given by the weight and volume capacities and the number of requests. For the directly aggregated input structure,
the features are derived from the capacities, number of items, and aggregate statistics of the items, namely, the weights, volumes, and off-loading costs are used.  For each of these, the min, max, mean, median, standard deviation, and 1st/3rd quartiles are computed.

\subsubsection{RL Training}

For RL, we have two distinct state space representations.  For DQN-L, we consider a linear state space with the following features: the incoming request index and price ratio, the number of accepted requests of each request class, the expected remaining weight and volume capacities, and the period.  For DQN-S, we consider a permutation-invariant model similar to the previous section with request features, $x_{r}$, given by the item class, the item class, price ratio, expected weight, expected volume, expected cost, expected price, and expected chargeable weight.  For the carrier features, $x_{c}$, we consider the previously stated request features for the incoming request, the expected weight and volume capacities, the expected utilized weight and volume, and the period.  Similar to the DiL case, the period can either be encoded with an integer or a one-hot vector.  

One important note for the training of DQN is that for the ACM instances, the EoHP is uncertain until the end of the booking period.  This results in a reward with high-variance depending on the realization of weights and volumes and not the trajectory itself.  To mitigate this, we sample 50 realizations of the weight and volume then take the average over the approximate EoHP costs as predicted by the supervised learning model.  For details on the choice of 50 realizations see Appendix~\ref{app:N_FINAL_OBS}.  This highlights an important area where supervised learning provides a significant advantage, even for the easier EoHP, as computing the 50 end-of-horizon costs can be done with efficient parallel matrix operations, whereas other approaches would need to access significantly more computing resources in order to parallelize a heuristic or exact algorithm across the different realizations.  Without additional sampling of EoHP costs, training is less stable and achieves significantly worse results.

\section{Computational Results}
\label{section:results}
This section reports the results of our computational experiments.  It is structured as follows. First, Section~\ref{sec:RES:setup} provides details about the computational setup and instances. Then, Section~\ref{sec:SL}, presents the supervised learning results, that is, the first phase of our methodology.  Section~\ref{sec:controlAlgosBaselines} outlines the specific parameters of the control algorithms and the baselines. Finally, we report and discuss our main results in Section~\ref{sec:ControlResults}.

\subsection{Computational Setup and Instances}
\label{sec:RES:setup}

All experiments were run on a computing cluster equipped with an Intel Xeon E5-2683 CPU and an Nvidia Tesla P100 GPU operating with
64GB of RAM. Gurobi 9.1.2 \citep{gurobi} handles  MILP solving,
Scikit-learn 1.0.1 \citep{scikit-learn} and Pytorch 1.10.0 \citep{NEURIPS2019_9015} handles supervised learning.  We build our implementation of RL based on the DQN implementation in RLLib \citep{liang2018rllib}.

In our experiments, we include eight and four instances from the DiL and ACM applications, respectively.  Here, an instance refers to an MDP, i.e., \eqref{eq:exact_mdp}-\eqref{eq:exact_boundary}, and  we refer to a single realization of the uncertainty as an episode or trajectory. We analyze results for 1,000 trajectories of each problem instance.  

For DiL, we denote an instance as VRP\_X\_Y, where $\text{X}\in\{4,10,15,50\}$ is the number of locations and $\text{Y}\in\{\text{L}, \text{H}\}$ is the load factor defining the capacity in terms of the expected demand.  A value of L indicates a low load factor, which results in a higher capacity, while a value of H indicates the opposite.  For each instance, locations are generated uniformly at random.  The number of time periods is 20, 30, 50, and 100 for the instances with 4, 10, 15, and 50 locations, respectively.  We choose the profits and probability of requests such that high-value requests have a higher probability of occurring later in the booking period.  A more detailed description of all the instances is available in Appendix~\ref{app:inst:vrp}.

For ACM, we denote an instance as CM\_X\_Y, where $\text{X}\in\{0.5, 1.0\}$ is the expected weight capacity divided by the expected total weight of all requests received during the booking period and $\text{Y}\in\{0.5, 1.0\}$ is the expected volume capacity divided by the expected total volume of all requests received during the booking period.  These instances are all generated with exactly the same parameters as the single-leg experiments in  \cite{barz2016air}. They have 60 time periods and 24 item classes.  Similar to the DiL instances, the profits and probabilities of requests are such that high-value requests have a higher probability of occurring later in the booking period.  A more detailed description of all the instances is available in Appendix~\ref{app:inst:cm}.

\subsection{Supervised Learning Results}
\label{sec:SL}

This section presents the numerical results of the supervised learning component of our methodology.  Specifically, Table~\ref{tab:sl_metrics} reports the mean absolute error (MAE) on both the train and validation sets for our permutation-invariant neural model (NN) described in Section~\ref{sec:EOH} as well as the time dedicated to training models and generating data. We compare with two baselines that use linear features derived from aggregating statistics of the requests, namely linear regression (LR) and random forests (RF) \citep{breiman2001random}.   

To generate datasets, we simulate random policies as described in Section~\ref{sec:EOH} to obtain 5,000 training and 1,000 validation samples of end-of-horizon state-cost pairs.  This process is done for each instance.  The choice of dataset size is motivated by the results in Appendix~\ref{app:MD:DSS}, which show diminishing returns when generating samples beyond these cardinalities.  All data generation was done in parallel with 32 threads, which is reflected in the reported computing times. The last column of Table~\ref{tab:sl_metrics} indicates that the time to generate data is quite small for all problem instances, ranging between 49 and 785 seconds.  As the datasets are generated offline before policies are computed, we consider these figures to be small.

For LR and RF the Scikit-learn implementation of 10-fold cross-validation is used for hyperparameter selection.  For the NN, we perform random searches over 50 configurations for each instance.  Detailed results of the best configurations from random search and the variance between configurations are reported in Appendix~\ref{app:MD:SL}. In summary, the variance in performance is relatively small, so it does not appear that hyperparameter configuration is particularly needed in this case.

\begin{table*}[t]\centering\resizebox{0.9\textwidth}{!}{
\begin{tabular}{l|rrr|rrr|rrr|r}
\toprule
Instance & \multicolumn{3}{c}{Train MAE}  & \multicolumn{3}{c}{Validation MAE}  & \multicolumn{3}{c}{Times (seconds)} & \\
\cmidrule{2-11}
 & LR & RF & NN & LR\  & RF\  & NN\  & LR\ \  & RF\ \  & NN\ \  & Data \\
\midrule
VRP\_4\_L & 1.46 & 0.62 & \textbf{0.03} & 1.47 & 0.69 & \textbf{0.09} & \textbf{3.31} & 78.00 & 355.71 & 175.72 \\
VRP\_10\_L & 2.56 & 0.90 & \textbf{0.74} & 2.57 & 1.32 & \textbf{0.81} & \textbf{2.84} & 144.94 & 373.53 & 341.75 \\
VRP\_15\_L & 2.13 & 0.93 & \textbf{0.82} & 2.11 & 1.35 & \textbf{0.93} & \textbf{3.19} & 182.29 & 170.12 & 420.84 \\
VRP\_50\_L & 10.00 & \textbf{2.67} & 4.58 & 9.77 & 6.89 & \textbf{5.41} & \textbf{3.49} & 233.43 & 665.70 & 657.05 \\
VRP\_4\_H & 2.75 & 1.21 & \textbf{0.07} & 2.73 & 1.37 & \textbf{0.16} & \textbf{2.68} & 84.15 & 220.66 & 189.12 \\
VRP\_10\_H & 3.82 & 1.31 & \textbf{1.16} & 3.76 & 2.25 & \textbf{1.42} & \textbf{2.73} & 151.04 & 165.46 & 381.30 \\
VRP\_15\_H & 2.99 & \textbf{1.18} & 1.40 & 3.09 & 2.12 & \textbf{1.55} & \textbf{3.09} & 162.09 & 236.50 & 530.62 \\
VRP\_50\_H & 13.27 & \textbf{5.76} & 5.78 & 13.23 & 9.22 & \textbf{6.06} & \textbf{2.88} & 223.50 & 193.62 & 785.69 \\
\midrule
CM\_0.5\_0.5 & 2034.86 & 141.73 & \textbf{55.53} & 1948.74 & 286.98 & \textbf{66.60} & \textbf{3.71} & 355.87 & 338.95 & 50.10 \\
CM\_1.0\_0.5 & 2018.82 & 90.74 & \textbf{27.20} & 2070.80 & 192.70 & \textbf{37.03} & \textbf{4.74} & 353.21 & 334.58 & 49.68 \\
CM\_0.5\_1.0 & 2013.31 & 83.59 & \textbf{26.19} & 2051.96 & 182.99 & \textbf{25.68} & \textbf{3.07} & 289.57 & 337.01 & 49.95 \\
CM\_1.0\_1.0 & 1966.05 & 93.93 & \textbf{30.01} & 1939.49 & 168.90 & \textbf{32.64} & \textbf{4.23} & 369.74 & 199.81 & 56.89 \\
\bottomrule
\end{tabular}}
\caption{Validation MAE, training times, and data generation times.  }
\label{tab:sl_metrics}
\end{table*}

The results in Table~\ref{tab:sl_metrics} show that NN consistently outperforms the baselines in validation MAE, often by an order of magnitude.  Whereas this does come at the cost of increasing training times, we argue that, considering that the training is performed offline, the considerable improvement in performance justifies the use of NN. Hence, we use NN as the EoHP cost predictor for all instances.

\subsection{Control Algorithms and Baselines} \label{sec:controlAlgosBaselines}

For our control algorithms, we implement two variants of DQN: DQN-L and DQN-S as described in Section~\ref{sec:RL:motivation}.  The DQN control results reported in the next section are obtained using random search over 100 configurations and selecting the model that obtains the highest mean reward over 100 validation trajectories.  Detailed results of the model selection and variance are reported in Appendix~\ref{app:MD:DQN}.  In addition, we report the mean reward on the validation trajectories throughout the DQN training in Appendix~\ref{app:DQN_CURVES}.

When deploying the trained DQN on the evaluation instances we additionally impose the capacity restrictions by simply rejecting requests that will violate the constraints (in expectation for the cargo management case).  Whereas this is not strictly necessary, given that a violation of constraints results in a significant penalty, we did empirically notice a slight improvement when imposing constraints during evaluation.  This may be directly handled by more sophisticated RL algorithms that can explicitly impose hard constraints \citep{donti2021enforcing,Ryu2020CAQL:}.

In terms of baselines, for DiL, we compare with BLP and BLPR as described in Section~\ref{sec:VRP:SOTA}, and for ACM we compare with DLP and DPD as described in Section~\ref{sec:CM:SOTA}.  We note that there are other baselines for ACM, such as randomized linear programming.  However, from \cite{barz2016air}, DPD outperforms all other algorithms for every instance.  DLP is included as the implementation is required for DPD. 
In our implementation of BLP/BLPR, we use a branch-and-cut procedure for the subtour elimination constraints \eqref{eq:vrp_sec} to solve the problems to optimality quite efficiently rather than relying on a fixed number on incumbent solutions as it is done in \cite{ggm}. In addition, we implement a naive first-come-first-serve (FCFS) baseline, that simply accepts requests until a capacity limit is reached.

\subsection{Control Results} \label{sec:ControlResults}

We now turn our attention to the main results of the paper. We compare the quality and speed of the control policies computed with our methodology to the baselines.  In Tables~\ref{tab:vrp_control_results} and~\ref{tab:cm_control_results}, we report mean profits, as well as online and offline computing timesfor the DiL and ACM applications, respectively.  The profit and online computing times are averaged over 1,000 trajectories, while the offline time is incurred once before any policy evaluation is done.  Note that the offline computing time does not include the time to generate the data and train the supervised learning models.

\begin{figure}[h!]
    \centering
    \includegraphics[width=1.0\textwidth]{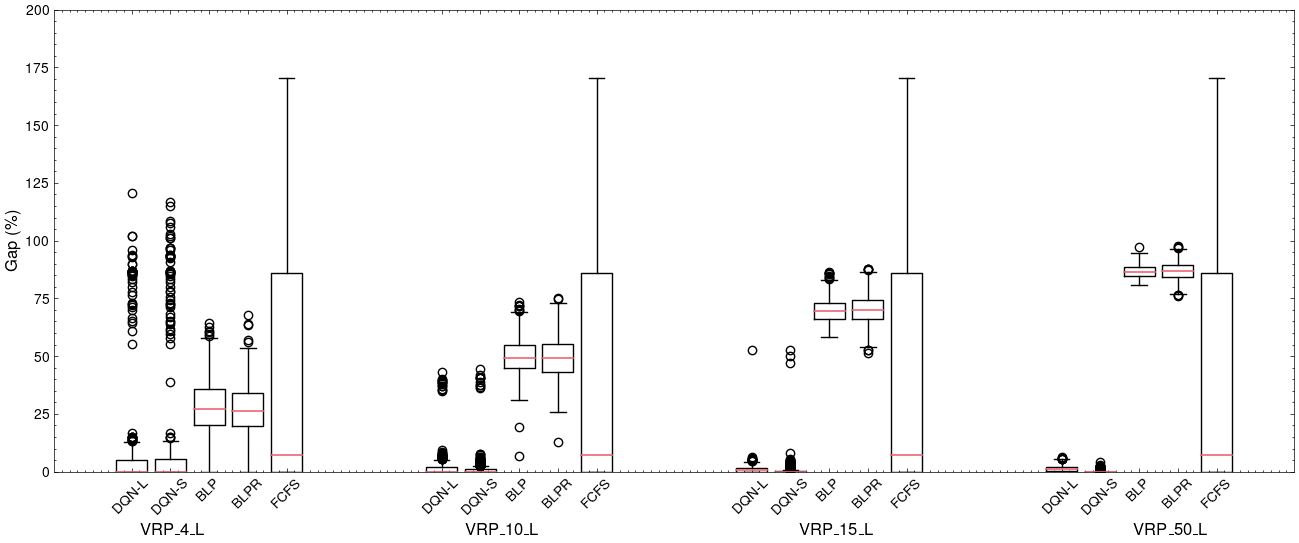}
    \caption[Gap to best known solution for DiL (VRP\_X\_L).]%
    {{\small Gap to best known solution for DiL (VRP\_X\_L). }}    
    \label{fig:vrp_L_boxplot}
\end{figure}

\begin{figure}[h!]
    \centering
    \includegraphics[width=1.0\textwidth]{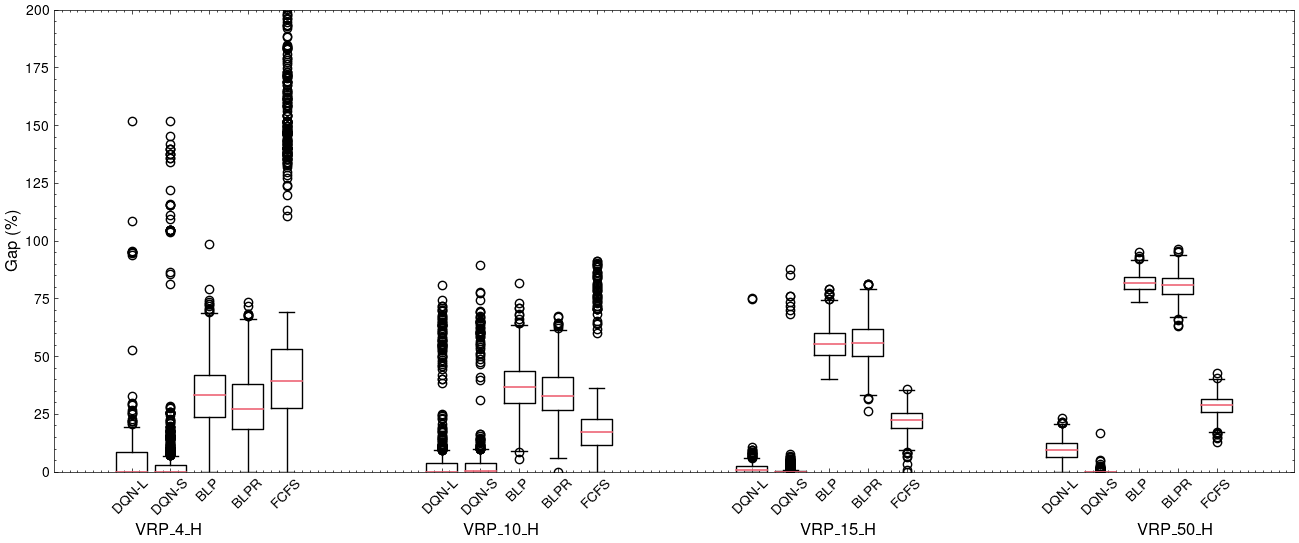}
    \caption[Gap to best known solution for DiL (VRP\_X\_H). ]%
    {{\small Gap to best known solution for DiL (VRP\_X\_H).}}    
    \label{fig:vrp_H_boxplot}
\end{figure}

\begin{figure}[h!]
    \centering
    \includegraphics[width=1.0\textwidth]{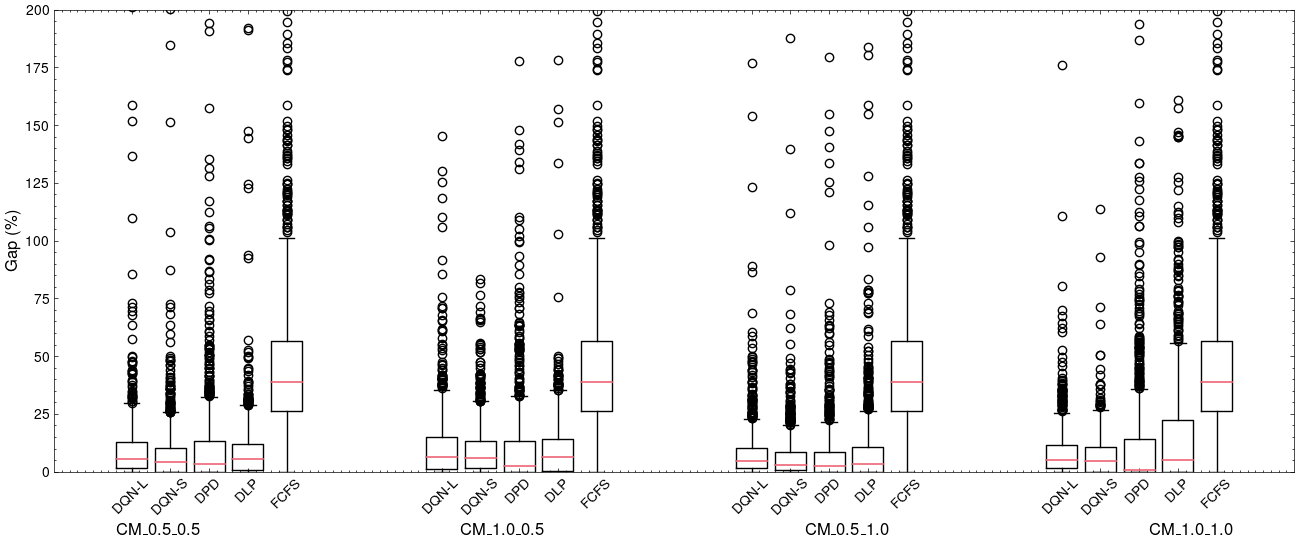}
    \caption[Gap to best known solution for ACM.]%
    {{\small Gap to best known solution for ACM.}}    
    \label{fig:cm_boxplot}
\end{figure}

In addition to the tables, we report trajectory-specific gaps to the best-known solutions as box plots in Figures~\ref{fig:vrp_L_boxplot}-\ref{fig:cm_boxplot}.  The gap is computed as  $ 100 \cdot \frac{\text{Best Reward} - \text{Algorithm Reward}}{\text{Best Reward}}$, where Best Reward is the best-known solution between all the algorithms for a given trajectory, and Algorithm Reward is the reward for a specific algorithm on the same trajectory.  Whereas this is a non-standard way of reporting the gaps, we use Best Reward as opposed to Algorithm Reward in the denominator in view of the significant difference in performance between our approach and the baselines for the DiL application.

\begin{table*}[t]\centering\resizebox{1.0\textwidth}{!}{
\begin{tabular}{l|rrrrr|rrrrr|rrrrr}
\toprule
Instance & \multicolumn{5}{c}{Mean Profits} & \multicolumn{5}{c}{Online Time} & \multicolumn{5}{c}{Offline Time} \\
\cmidrule{2-16}
 & DQN-L & DQN-S & BLP & BLPR & FCFS & DQN-L\  & DQN-S\  & BLP\  & BLPR\  & FCFS\  & DQN-L\ \  & DQN-S\ \  & BLP\ \  & BLPR\ \  & FCFS\ \  \\
\midrule
VRP\_4\_L & \textbf{102.41} & 101.17 & 77.12 & 78.68 & 70.23 & 0.82 & 0.85 & \textbf{0.01} & 0.12 & 0.87 & 2178.41 & 3118.10 & 0.09 & \textbf{0.08} & - \\
VRP\_10\_L & 246.70 & \textbf{247.99} & 125.57 & 127.69 & 243.45 & 2.31 & 2.73 & \textbf{0.08} & 0.23 & 2.22 & 3730.99 & 4361.60 & \textbf{0.12} & 0.14 & - \\
VRP\_15\_L & 465.31 & \textbf{468.01} & 142.86 & 141.39 & 451.04 & 2.75 & 2.80 & \textbf{0.19} & 0.33 & 2.64 & 6143.76 & 8355.47 & 0.21 & \textbf{0.20} & - \\
VRP\_50\_L & 1409.08 & \textbf{1426.69} & 190.92 & 186.30 & 1372.16 & 3.99 & 3.49 & \textbf{0.88} & 2.10 & 3.25 & 12565.44 & 19482.41 & 1.15 & \textbf{1.10} & - \\
VRP\_4\_H & \textbf{77.13} & 75.30 & 52.72 & 57.85 & 32.44 & 0.75 & 0.82 & \textbf{0.01} & 0.14 & 0.87 & 2171.51 & 2655.38 & 0.31 & \textbf{0.24} & - \\
VRP\_10\_H & 152.69 & \textbf{154.08} & 102.15 & 107.05 & 129.75 & 2.00 & 2.10 & \textbf{0.11} & 0.46 & 2.01 & 5044.56 & 4257.89 & 0.82 & \textbf{0.51} & - \\
VRP\_15\_H & 314.55 & \textbf{316.15} & 142.86 & 141.39 & 248.36 & 2.37 & 2.34 & \textbf{0.18} & 0.33 & 2.40 & 5739.93 & 6890.85 & 0.21 & \textbf{0.20} & - \\
VRP\_50\_H & 823.07 & \textbf{909.24} & 165.91 & 178.51 & 647.29 & 2.78 & 2.43 & \textbf{0.95} & 2.17 & 3.15 & 13018.87 & 16230.71 & 1.38 & \textbf{1.12} & - \\
\bottomrule
\end{tabular}}
\caption{DiL control results.  Profit and online computing time averaged over 1,000 trajectories.  All times in seconds.}
\label{tab:vrp_control_results}
\end{table*}

\begin{table*}[t]\centering\resizebox{1.0\textwidth}{!}{
\begin{tabular}{l|rrrrr|rrrrr|rrrrr}
\toprule
Instance & \multicolumn{5}{c}{Mean Profits} & \multicolumn{5}{c}{Online Time} & \multicolumn{5}{c}{Offline Time} \\
\cmidrule{2-16}
 & DQN-L & DQN-S & DPD & DLP & FCFS & DQN-L\  & DQN-S\  & DPD\  & DLP\  & FCFS\  & DQN-L\ \  & DQN-S\ \  & DPD\ \  & DLP\ \  & FCFS\ \  \\
\midrule
CM\_0.5\_0.5 & 6094.05 & \textbf{6136.34} & 6111.81 & 6100.18 & 3816.06 & 0.21 & 1.33 & 0.03 & \textbf{0.03} & 0.04 & 8181.33 & 8360.47 & 273.48 & 0.05 & \textbf{0.00} \\
CM\_1.0\_0.5 & 6417.52 & 6439.14 & \textbf{6559.94} & 6433.58 & 4341.09 & 0.20 & 1.30 & 0.03 & \textbf{0.02} & 0.03 & 6361.24 & 8596.11 & 431.28 & 0.07 & \textbf{0.00} \\
CM\_0.5\_1.0 & 6430.19 & 6530.65 & \textbf{6538.72} & 6421.70 & 4282.58 & 0.18 & 0.99 & 0.03 & \textbf{0.03} & 0.03 & 5765.68 & 8861.07 & 363.32 & 0.06 & \textbf{0.00} \\
CM\_1.0\_1.0 & 9913.02 & \textbf{9967.68} & 9851.77 & 9306.05 & 9306.05 & 0.16 & 1.02 & \textbf{0.04} & 0.04 & 0.04 & 7127.56 & 8760.90 & 528.42 & 0.06 & \textbf{0.00} \\
\bottomrule
\end{tabular}}
\caption{ACM control results.  Profit and online computing time averaged over 1,000 trajectories.  All times in seconds.}
\label{tab:cm_control_results}
\end{table*}

Based on these results, a number of findings emerge. For the DiL application, DQN-L and DQN-S perform significantly better than BLP and BLPR in both the low and high-load factor settings.  The primary reason for this is that booking-limit policies from BLP and BLPR are simply not adequate in the case where locations are not expected to receive a large number of requests.  In fact, BLP and BLPR are significantly worse than FCFS in this case, since for the majority of the locations the booking-limit policy provides a threshold of zero.  

DQN is able to achieve a control of better quality than FCFS that improves with increased load factor.  This is in line with our expectations as the value of making strategic decisions increases as less requests can be accepted.  In terms of our policies, DQN-S generally outperforms DQN-L, especially as the number of locations increases.  We hypothesize that this is due to DQN-S having a better representation of the locations and features of the previously accepted requests, especially as the number of locations increases. 

For the ACM application, we can see that DQN-L and DQN-S are able to achieve policies that are quite close to state of the art. 
Specifically, when compared to DPD, the percentage improvement in mean profit ranges between -2.2\% and 1.1\% (i.e., DQN provides a slight improvement compared to DPD).
Since DPD is a very strong baseline that leverages problem-specific structure, we did not expect to outperform it. The fact that our general methodology achieves a comparable performance is evidence of the potential of RL for booking control.
Moreover, compared to DLP, which is less problem-specific than DPD, we can see that we are able to obtain higher-quality policies in three out of the four instances for both DQN-L and DQN-S.  
Comparing DQN-L to DQN-S, we can see that there is no significant difference between the two for this instance size.

As expected, the offline computing time is much larger for DQN-L and DQN-S than the other algorithms. Nevertheless, they are still reasonable for an end-to-end RL approach as they take less than six hours in the worst case when combining data generation, supervised learning training, and DQN training. The algorithm can hence be retrained in an overnight computing time budget. 
We close this section by discussing the reduction in computing time resulting from the supervised learning phase of our approach. 
If the supervised learning prediction were unavailable, solving or approximating 45,000 and 2,250,000 EoHPs would be required to train and evaluate a single DQN configuration for the DiL and ACM applications, respectively.  The larger number -- 2,250,000 EoHPs -- is required as we sample 50 EoHP problems for each trajectory in the ACM application.

Let us analyze the time spent evaluating the EoHP within a DQN training run. We compare them in a completely sequential manner. We take VRP\_50\_H as an example. Solving a single EoHP takes 
5.03 seconds on average with FILO whereas our ML approximation with bin-packing takes 
0.02 seconds.
We extrapolate to compute the total time by multiplying these figures with the number of EoHPs evaluated during DQN training (45,000), which results in 226278.72 seconds ($\sim 62.8$ hours) for FILO and 
1093.57 seconds for our approximation.  

Similarly for CM\_1.0\_1.0,  solving a single EoHP takes 
0.31 seconds on average while our ML approximation takes 
0.5 milliseconds.  Multiplying these figures with the number of EoHPs evaluated during DQN training (2,250,000) results in 819216.00 seconds ($\sim 227.5$ hours) for the exact algorithm and 
1129.57 seconds for our approximation. We note that in our actual implementation this cost is roughly 50 times lower, as all 50 EoHP realizations are fed through the neural network in a single batch.  

Overall, the time spent approximating the EoHP cost during simulation with ML is several orders of magnitude smaller than with the algorithms used to compute the labels.  It is very unlikely that more efficient implementations or approximations of the EoHPs would achieve a time close to the nearly negligible amount of time required by our ML approximation. Importantly, this time is expected to be nearly invariant with respect to the specific EoHP structure. This hence opens up the possibility to use less crude approximations of the end-of-horizon operational decision-making problem. In turn, this may result in policies that better reflect real operations (e.g., a three-dimensional knapsack instead of a vector packing problem formulation in the ACM EoHP). The latter is not evaluated in the performance metrics reported here as it would require a field study.

\section{Conclusion}
\label{section:conclusion}

In this work, we focused our attention on freight booking control and proposed a two-phase approach, that first estimates the cost of the EoHP via supervised learning, then deploys the trained model to significantly speed up simulation for RL algorithms.  Our approach can be extended to any freight booking control problem, regardless of the EoHP, as long as one can adequately compute an exact or relevant approximation to the objective function value.  In addition, since a trained model can provide fast estimates for any EoHP, our approach can scale to problems with hard EoHPs.  Through our computational experiments, we demonstrated that our approach is able to achieve control policies that surpass the state-of-the-art ones by a significant margin in a DiL application and a policy comparable with a strong problem-specific algorithm in ACM.  

In terms of future work,
we identify five potential directions.  
First, as this work uses DQN for computing control policies, the extension to other RL algorithms may lead to even better performing control.  In addition, with the capacity constraint structure of the problem, one can consider RL algorithms that can explicitly enforce these constraints such as in \cite{donti2021enforcing} and \cite{Ryu2020CAQL:}.  

Second, determining better baselines or bounds \citep[e.g.,][]{brown2010information, brown2022information} would be useful in understanding the performance of our approach in regard to the application to DiL since the policies computed by BLP and BLPR were not of high enough quality to provide an informative comparison.

Fourth, it may be worth exploring an iterative improvement of the supervised learning model during the training of the DQN.  As the DQN trajectories improve with training, having a representative sample of those EoHPs may further improve the EoHP predictor.  We note that this process could be performed in parallel with the training of the DQN, so that the computing cost associated with these improvements would be negligible.

Last, the broadest direction for future work is the extension of our methodology to other freight booking control problems and more generally other sequential decision-making problems that involve solving discrete optimization problems during simulation.

\section*{Acknowledgment}
We are grateful to Luca Accorsi and Daniele Vigo who provided us access to their VRP code FILO \citep{accorsi2021fast}, which has been very useful for our research. We thank Francesca Guerriero for sharing details on the implementation of the approach in \cite{ggm}. We would also like to thank Pierre-Luc Bacon for valuable comments on this work, and we express our gratitude to Eric Larsen and Rahul Patel for their comments that helped us improve the manuscript. Computations were made on the clusters managed by Calcul Qu\'ebec and Compute Canada.

\clearpage

\appendix
\label{app:app}

\section{Acronyms \& Symbols}

List of symbols used in the paper with their brief description. 

\begin{table*}[h]\centering\resizebox{0.4\textwidth}{!}{
\begin{tabular}{c|l}
\toprule
    ACM      &  Airline cargo management \\
    ADP      &  Approximate dynamic programming  \\
    BLP      &  Booking limit policy  \\
    BLPR     &  Booking limit policy with reoptimization  \\
    DiL      &  Distributional Logistics \\
    DLP      &  Deterministic linear programming  \\
    DPD      &  Dynamic programming decomposition  \\
    DQN      &  Deep Q-learning  \\
    DQN-L    &  DQN with value function approximation from linear state  \\
    DQN-S    &  DQN with value function approximation from set-based state  \\
    EoHP     &  End of horizon problem  \\
    FCFS     &  First-come-first-serve  \\
    FILO     &  Fast Iterated local search Localized Optimization (VRP solver)  \\
    LP       &  Linear programming  \\
    LR       &  Linear regression  \\
    MAE      &  Mean absolute error  \\
    MDP      &  Markov decision process  \\
    MILP     &  Mixed-integer linear programming  \\
    ML       &  Machine learning  \\
    MTP      &  Bin packing solver  \\
    NN       &  Neural Network (specifically the permutation-invariant model)  \\
    RF       &  Random forests  \\
    RL       &  Reinforcement learning  \\
    RM       &  Revenue management  \\
    TSP      &  Travelling salesperson problem  \\
    VRP      &  Vehicle routing problem  \\
\bottomrule
\end{tabular}}
\caption{List of acronyms}
\label{tab:acronyms}
\end{table*}

\begin{table*}[h]\centering\resizebox{0.45\textwidth}{!}{
\begin{tabular}{c|l}
\toprule
\multicolumn{2}{c}{Booking Control}  \\
\midrule
    $\mathcal{N}$       & Set of requests   \\
    $n$                 & Cardinality the set of requests   \\
    $j$                 & Index for the set of requests   \\
    $T$                 & Number of periods for requests   \\
    $t$                 & Index for the periods   \\
    $p_j^t$             & Probability of request $j$ at time $t$  \\
    $p_0^t$             & Probability of no request at time $t$  \\
    $s_j^t$             & Number of accepted requests of class $j$ at time $t$  \\
    $\vect{s}^t$        & The state at time $t$, i.e. the vector of accepted requests   \\  
    $\mathcal{S}^t$     & Set of all possible states at time $t$  \\  
    $a_j$               & Action taken for a given request $j$   \\  
    $\vect{e}_j$        & Vector of all 0's and 1 at index $j$ \\
    $r_j$               & Revenue of request class $j$      \\
    $R(\cdot, \cdot)$   & Reward function ($R(j, a_j)=r_ja_j$)     \\
    $V_t(\cdot)$        & Value function\\
    $\Gamma(\cdot)$     & Cost of the EoHP \\
    $A, b$              & General constraints for DLP  \\
    $\vect\mu$          & Vector of expected number of requests for DLP  \\
\midrule
\multicolumn{2}{c}{Approximation \& Machine Learning}  \\
\midrule
    $\phi$              & General approximation for EoHP  \\ 
    $g$                 & Feature mapping for EoHP state   \\ 
    $\tilde{V}_t(\cdot)$ & Approximate value function  \\ 
    $x_r$               & Features for a given request \\
    $x_c$               & Features for a carrier \\
    $\Psi^1$            & Request embedding network \\
    $\Psi^2$            & Request decoding network \\
    $\Phi$              & Prediction network for EoHP \\
\midrule
\multicolumn{2}{c}{Distributional Logistics}  \\
\midrule
    $Q$                 & Capacity of each vehicle  \\
    $K_0$               & Number of no-cost vehicles \\
    $K$                 & Number of utilized vehicles \\
    $\mathcal{V}$       & Set of locations, $\mathcal{N}\cup\{0\}$\\
    $\mathcal{A}$       & Set of arcs \\
    $c_{ij}$            & Cost of an arc \\
    $z^*(\cdot, \cdot)$ & EoHP optimal objective \\
    $C$                 & Cost of additional vehicles for $K>K_0$  \\
    $\alpha_{ij}^k$     & Binary variable for vehicle $k$ using arc $(i,j)$  \\
    $\beta_{i}^k$       & Binary variable for vehicle $k$ visiting location $i$  \\
    $q_{j}^k$           & Variable for the quantity of accepted requests at location $j$ on vehicle $k$        \\
    $x_j$               & DLP variable for the number of accepted requests for location $j$ \\
\midrule
\multicolumn{2}{c}{Airline Cargo Management}  \\
\midrule
    $Q_w$               & Realized weight capacity \\
    $Q_v$               & Realized volume capacity \\
    $w_{ji}$            & Realized weight of the $i$-th accepted request of class $j$ \\
    $v_{ji}$            & Realized volume of the $i$-th accepted request of class $j$ \\
    $\eta$              & Weight-volume ratio   \\
    $c_j$               & Offloading cost of request class $j$    \\
    $z^*(\cdot)$        & EoHP optimal objective \\
    $E[Q_w]$            & Expected weight capacity \\
    $E[Q_v]$            & Expected volume capacity \\
    $E[w_j]$            & Expected weight of request class $j$ \\
    $E[v_j]$            & Expected volume of request class $j$ \\
    $y_{ji}$            & Binary variable for loading $i$-th request of class $j$ for EoHP \\
    $x_j$               & DLP variable for the number of accepted for requests class $j$ \\
    $y_j$               & DLP variable for the loaded quantity of requests class $j$ \\
    $z_j$               & DLP variable for the offloaded quantity of requests class $j$ \\
\bottomrule
\end{tabular}}
\caption{List of symbols}
\label{tab:symbols}
\end{table*}

\section{Description of Instances}
\label{app:inst}
In this section, we provide detailed descriptions of how the instances were generated for both the DiL and ACM applications.

\subsection{Distributional Logistics}
\label{app:inst:vrp}

We detail the parameters that define the 8 sets of instances in the DiL application.  In each of the instances the capacity is determined to be inversely proportional to the demand outstanding at the locations.  Specifically, a load factor, $LF$, is used, and the capacity is given by

\begin{equation}
    Q = \floor[\Big]{ \frac{\sum_{j\in\mathcal{N}\cup\{0\}} \sum_{t=1}^T p_j^t}{K_0\cdot LF} }.
\end{equation}

This load factor is used to determine capacity differences between the instances ending in L and H, where L corresponds to a low load factor and H corresponds to a high load factor.  This is done to create two distinct instances where accept/reject decisions become more crucial as the load factor is increased.  The locations are determined by sampling 2D coordinates uniformly at random in a box.  For each instance, we specify the number of locations, periods, load factors, number of no cost vehicles, additional vehicle costs, and location bounds in Table~\ref{tab:vrp_inst_params}.

\begin{table*}[h!]\centering\resizebox{0.7\textwidth}{!}{
\begin{tabular}{l|cccccc}
\toprule
{Instance}      & {\# locations}    &  {\# periods}   & {Load factor} & {\# No cost vehicles}  & {Additional vehicle cost }  & {Location bounds}                 \\
\midrule
VRP\_4\_L       &  4    &  20   &  1.1  &  2    &  100   & $[0,10]$  \\                
VRP\_10\_L      &  10   &  30   &  1.1  &  3    &  100   & $[0,10]$  \\                
VRP\_15\_L      &  15   &  50   &  1.1  &  3    &  250   & $[0,10]$  \\             
VRP\_50\_L      &  50   &  100  &  1.1  &  3    &  600   & $[0,50]$  \\            
VRP\_4\_H       &  4    &  20   &  1.8  &  2    &  100   & $[0,10]$  \\         
VRP\_10\_H      &  10   &  30   &  1.8  &  3    &  100   & $[0,10]$  \\            
VRP\_15\_H      &  15   &  50   &  1.8  &  3    &  250   & $[0,10]$  \\           
VRP\_50\_H      &  50   &  100  &  1.8  &  3    &  600   & $[0,50]$  \\             
\bottomrule
\end{tabular}}
\caption{Summary of instances generating parameters for DiL.}
\label{tab:vrp_inst_params}
\end{table*}

For the profit we aim to generate instances which have a higher probability of high profit requests later in the booking period.  We now detail how the probability and pricing parameters are defined to generate each instance in the DiL application.

\textbf{VRP\_4\_\{L,H\}.}
The profits for accepting a request from each location are defined by $r_1=4$, $r_2=8$, $r_3=12$, $r_4=16$.  The probability of no request is $p_0^t = 0.10$, for $t\in\{1,\ldots,T\}$.  The initial probability of a request for each location is $p_1^1 = 0.45$, $p_2^1 = 0.40$, $p_3^1 = 0.10$, $p_4^1 = 0.05$.  For the remainder of the periods, the probability of a request for each location is given by $p_j^{t+1} = p_j^{t} - 0.01$ for $j \in\{ 1,2\}$ and $p_j^{t+1} = p_j^{t} + 0.01$ for $j \in\{3,4\}$.

\textbf{VRP\_10\_\{L,H\}.}
The profits for accepting a request from each location are defined by $r_j=10$ for $j\in\{1,\ldots,4\}$, $r_j=12$ for $j\in\{5,\ldots,8\}$, and $r_j=20$ for $j\in\{9,10\}$.  The probability of no request is $p_0^t = 0.10$, for $t\in\{1,\ldots,T\}$.  The initial probability of a request for each location is $p_j^1 = 0.125$ for $j\in\{1,\ldots 4\}$, $p_j^1 = 0.075$ for $j\in\{5,\ldots,8\}$, and $p_j^1 = 0.05$ for $j\in\{9,10\}$.  For the remainder of the periods, the probability of a request for each location is given by $p_j^{t+1} = p_j^{t} - 0.001$ for $j\in\{1,\ldots4\}$,  $p_j^{t+1} = p_j^{t}$ for $j\in\{5,\ldots,8\}$, and $p_j^{t+1} = p_j^{t} + 0.002$ for $j\in\{9,10\}$. 

\textbf{VRP\_15\_\{L,H\}.}
The profits for accepting a request from each location are defined by $r_j=10$ for $j\in\{1,\ldots,5\}$, $r_j=12$ for $j\in\{6,\ldots,10\}$, and $r_j=20$ for $j\in\{11,\ldots,15\}$.  The probability of no request is $p_0^t = 0.10$, for $t\in\{1,\ldots,T\}$.  The initial probability of a request for each location is $p_j^1 = 0.10$ for $j\in\{1,\ldots,5\}$, $p_j^1 = 0.06$ for $j\in\{6,\ldots,10\}$, and $p_j^1 = 0.02$ for $j\in\{11,\ldots,15\}$.  For the remainder of the periods, the probability of a request for each location is given by $p_j^{t+1} = p_j^{t} - 0.001$ for $j\in\{1,\ldots,5\}$,  $p_j^{t+1} = p_j^{t}$ for $j\in\{6,\ldots,10\}$, and $p_j^{t+1} = p_j^{t} + 0.001$ for $j\in\{11,\ldots,15\}$. 

\textbf{VRP\_50\_\{L,H\}.}
The profits for accepting a request from each location are defined by $r_j=15$ for $j\in\{1,\ldots,30\}$, $r_j=22$ for $j\in\{31,\ldots,40\}$, and $r_j=30$ for $j\in\{41,\ldots,50\}$.  The probability of no request is $\lambda_0^t = 0.10$, for $t\in\{1,\ldots,T\}$.  The initial probability of a request for each location is $p_j^1 = 0.0166$ for $j\in\{1,\ldots,30\}$, $p_j^1 = 0.03$ for $j\in\{31,\ldots,40\}$, and $p_j^1 = 0.01$ for $j\in\{41,\ldots,50\}$.  For the remainder of the periods, the probability of a request at for each location is given by $p_j^{t+1} = p_j^{t} - 0.0001$ for $j\in\{1,\ldots,30\}$,  $p_j^{t+1} = p_j^{t}$ for $j\in\{31,\ldots,40\}$, and $p_j^{t+1} = p_j^{t} + 0.0003$ for $j\in\{41,\ldots,50\}$.

\subsection{Airline Cargo Management}
\label{app:inst:cm}

We describe the parameters  defining the ACM instances in the second application.  These instances are taken directly from \cite{barz2016air}.  
In our setting the number of periods is given by $T=60$.  The parameter which determines the chargeable weight is given by $\eta=0.6$.  The number of shipment classes is $n=24$.  As mentioned previously, whereas the realizations of the weight and volume are given at the end of the booking period, the expectations of each class are known and provided in Table~\ref{tab:air_expected}, where $\bar w_j$ and $\bar v_j$ denote the expected weight and volume of class $j$.  The realized weight and volume of a request of class $j$ at the end of a booking period are given by sampling a multivariate normal distribution with mean $\vect \mu = [\bar w_j, \bar v_j]$ and covariance matrix

\begin{align}
\Sigma &= \begin{pmatrix}
    (\sigma_j^w)^2 & 0.8 \sigma_j^w \sigma_j^v \\
     0.8 \sigma_j^w \sigma_j^v & (\sigma_j^v)^2 \\
\end{pmatrix},
\end{align}

where $\sigma_j^w = 0.25 \bar w_j$ and $\sigma_j^v = 0.25\bar v_j$.  

\begin{table*}[h!]\centering\resizebox{\textwidth}{!}{
\begin{tabular}{l|cccccccccccccccccccccccc}
\toprule
Class                           & 1  & 2 & 3 & 4 & 5 & 6 & 7 & 8 & 9 & 10 & 11  & 12 & 13 & 14 & 15 & 16 & 17 & 18 & 19 & 20 & 21 & 22 & 23 & 24   \\
\midrule
$\bar w_j$ (kg)                 & 50 & 50 & 50 & 50 & 100 & 100 & 100 & 100 & 200 & 200 & 250 & 240 & 300 & 400 & 500 & 500 & 1000 & 1500 & 2500 & 3500 & 70 & 70 & 210 & 210  \\ 
$\bar v_j$ ($1e^{4}$ cm$^3$)    & 30 & 29 & 27 & 25 & 59 & 58 & 55 & 52 & 125 & 119 & 100 & 147 & 138 & 179 & 235 & 277 & 598 & 898 & 1488 & 2083 & 244 & 17 & 700 & 52 \\
\bottomrule
\end{tabular}}
\caption{ACM -- Expected weights and volumes.}
\label{tab:air_expected}
\end{table*}

Let $w_j$ and $v_j$ be the realized weight and volume of a request of class $j$.  
The off-loading cost of the request is given by $c_j\cdot \max\{w_{j}, v_{j} / \eta\}$, where $c_j = 2.4$ for $j=1,\ldots, 24$. 
The profit is given by $\alpha \cdot \max\{w_{j}, v_{j} / \eta\}$, where $\alpha \in \{0.7, 1.0, 1.4\}$.  
Here, $\alpha$ is used in order to ensure that more profitable requests occur later in the booking period.  Table~\ref{tab:air_alpha} gives the probability that $\alpha$ takes each value over the time horizon as well as the probability that no request occurs.  Given that a request occurs, the probability of the request for each class is provided in Table~\ref{tab:air_class_probs}.

\begin{table*}[h!]\centering\resizebox{0.4\textwidth}{!}{
\begin{tabular}{l|ccc}
\toprule
Request                     &                & Periods Range    &      \\
\cmidrule{2-4}
                            & 1-20          &  21-40            & 41-60     \\
\midrule
$\text{Pr}(\alpha = 0.7)$          & 0.7           & 0.4               & 0.0       \\
$\text{Pr}(\alpha = 1.0)$          & 0.2           & 0.2               & 0.2       \\
$\text{Pr}(\alpha = 1.4)$          & 0.0           & 0.4               & 0.7       \\
$\text{Pr}(\text{No request})$      & 0.1           & 0.0               & 0.1       \\
\bottomrule
\end{tabular}}
\caption{ACM -- Profit per chargeable weight probabilities.}
\label{tab:air_alpha}
\end{table*}

\begin{table*}[h!]\centering\resizebox{\textwidth}{!}{
\begin{tabular}{l|cccccccccccccccccccccccc}
\toprule
Class                           & 1  & 2 & 3 & 4 & 5 & 6 & 7 & 8 & 9 & 10 & 11  & 12 & 13 & 14 & 15 & 16 & 17 & 18 & 19 & 20 & 21 & 22 & 23 & 24   \\
\midrule
Probability                 &  0.072 &  0.072 &  0.072 &  0.072 &  0.072 &  0.072 &  0.072 &  0.072 &  0.072 &  0.072 &  0.040 &  0.040 &  0.040 &  0.040 &  0.040 &  0.040 & 0.009 & 0.009 & 0.009 & 0.009 & 0.001 & 0.001 & 0.001 & 0.001  \\ 
\bottomrule
\end{tabular}}
\caption{ACM --  Class request probabilities.}
\label{tab:air_class_probs}
\end{table*}

Similar to the weight and volume of a request, the capacities are realized at the end of the booking period.  The expected weight and volume capacities are given as a ratio of the total expected weight and volume of the requests.  In the context of this work, we consider all single-leg instances from \cite{barz2016air}.  The expected weight capacity and volume capacity are denoted by $\bar Q_w$ and $\bar Q_v$, respectively.  At the end of the booking period, the realized weight and volume capacities are given by sampling a multivariate normal distribution with mean $\vect \mu = [\bar Q_w, \bar Q_v]$ and covariance matrix

\begin{align}
\Sigma &= \begin{pmatrix}
    (\sigma^{Q_w})^2 & 0.8 \sigma^{Q_w} \sigma^{Q_v} \\
     0.8 \sigma^{Q_w} \sigma^{Q_v} & (\sigma^{Q_v})^2 \\
\end{pmatrix},
\end{align}
where $\sigma^{Q_w} = 0.25 \bar Q_w$ and $\sigma^{Q_v} = 0.25\bar Q_v$.  

Finally, we detail the differences between the instances.  For instance CM\_X\_Y, X specifies the ratio of expected weight capacity in comparison to the expected weight of all requests, while Y specifies the same for volume.  For example, in instance CM\_0.5\_1.0, the expected total weight of all requests is twice as large as the expected weight capacity, and the expected total volume of all requests is equal to the expected volume capacity.

\section{Model Selection and Dataset Sizing}
\label{app:MD}

In this section, we report on the effects of variations in dataset sizes and hyperparameter configurations for both the supervised learning model for predicting the EoHP cost and the DQN models for control.

\subsection{Dataset Sizing}
\label{app:MD:DSS}

\begin{figure}[h!]
    \centering
    \includegraphics[width=0.5\textwidth]{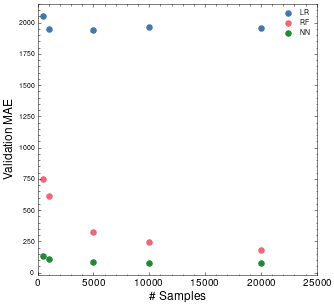}
    \caption[Dataset sizing in ACM (CM\_0.5\_0.5).]%
    {{\small Dataset sizing in ACM (CM\_0.5\_0.5).}}    
    \label{fig:cm_dss}
\end{figure}

As supervised learning is a core component of our methodology, determining an estimate of a good number of samples is important. Ideally this would be done for every instance under consideration.  However, this would impose a large computational burden, so we restrict our attention to a single instance setting, namely, CM\_0.5\_0.5 and generate training datasets with 500, 1,000, 5,000, 10,000, and 20,000 samples.  We then train LR, RF, and NN on each training set size and evaluate every model on the same 5,000 validation samples.  The validation MAEs for each model type and dataset size are reported in Figure~\ref{fig:cm_dss}.  From this figure, we can see that with LR and NN, there is no significant improvement beyond 5,000 samples, so for this reason we use 5,000 training samples when generating data for all of the instances.

\subsection{Supervised Learning Model Selection}
\label{app:MD:SL}

For the supervised learning models, LR and RF are tuned with random search and 10-fold cross validation using the implementation in Scikit-learn.  Only the best model is reported in the main paper.  For NN, random search is used for model selection.  The batch size, learning rate, activation function, and \# of epochs are all relatively standard hyperparameters, so we omit their description and only report their values.  However, the remainder are specific to our permutation-invariant model.  To simplify the hyperparameter search for the NN models, we consider a fixed-depth network and allow the width to be a configurable hyperparameter. 
Embed hidden dimension 1 and Embed dimension 2 indicate the width of the two layers of $\Psi_1$, Embed hidden dimension 3 is the width of $\Psi_2$ and  Concat is the width of the single layer in $\Phi$. $\Psi_1$, $\Psi_2$ and $\Phi$ are illustrated in Figure~\ref{fig:nn_arch}.

\begin{table*}[h!]\centering\resizebox{0.4\textwidth}{!}{
\begin{tabular}{l|c}
\toprule
{Hyper parameter name}                & {Parameter space}                       \\
\midrule
    Batch size                  & $\{16, 32, 64, 128\}$         \\
    Learning rate               & $[1e^{-5},1e^{-2}]$           \\
    Activation Function         & \{ReLU, Leaky Relu\}          \\
    \# Epochs                   & $1000$                        \\
    Embed hidden dimension 1    & $\{64, 128, 256\}$            \\
    Embed hidden dimension 2    & $\{64, 128, 256\}$            \\
    Embed hidden dimension 3    & $\{64, 128, 256\}$            \\
    Concat hidden dimension     & $\{256, 512, 1024\}$          \\
\bottomrule
\end{tabular}}
\caption{Random search hyperparameter space for the neural network used in supervised learning.  For values in \{\}, we sample with equal probability for each discrete choice.  For values in [], we sample a uniform distribution within the given bounds.  For single values, we keep it fixed across all configurations.}
\label{tab:random_search_space_sl}
\end{table*}

\begin{figure}[h!]
    \centering
    \includegraphics[width=1.0\textwidth]{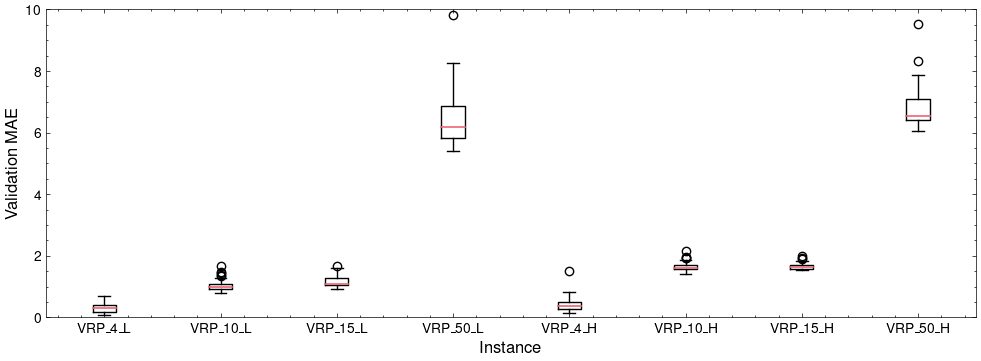}
    \caption[DiL box plots of MAE over random search configurations. ]%
    {{\small DiL box plots of MAE over random search configurations.}}    
    \label{fig:vrp_sl_config}
\end{figure}

\begin{figure}[h!]
    \centering
    \includegraphics[width=0.5\textwidth]{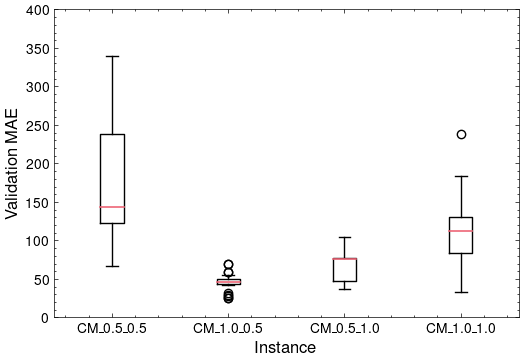}
    \caption[ACM box plots of MAE over random search configurations. ]%
    {{\small ACM box plots of MAE over random search configurations. }}    
    \label{fig:cm_sl_config}
\end{figure}

Within the search space specified in Table~\ref{tab:random_search_space_sl}, we perform random search over 50 configurations.  The best hyperparameters for each instance are reported in Table~\ref{tab:sl_best_config}.  In addition, we provide box plots showing the validation MAE over hyperparameter configurations in Figures~\ref{fig:vrp_sl_config} and \ref{fig:cm_sl_config} for the DiL and ACM problems, respectively. 
Each box plot shows the variations in MAE over the hyperparameter configurations for a particular instance. We conclude that we are generally able to achieve trained models of good quality reliably for the majority of the instances, 
with the exception of the larger DiL instances with 50 nodes.  It is worth noting here that even the worst performing models in the DiL case actually achieve a very low MAE, especially when compared to the target distribution, which is in the order $10^2$.  

\begin{table*}[t]\centering\resizebox{1.0\textwidth}{!}{
\begin{tabular}{l|rrrrrrr}
\toprule
{} & Learning rate & Activation function & Batch size & Embed hidden dimension 1 & Embed hidden dimension 2 & Embed hidden dimension 3 & Concat hidden dimension \\
\midrule
CM\_0.5\_0.5 &       0.00018 &          leaky\_relu &         32 &                      256 &                      128 &                      128 &                     512 \\
CM\_1.0\_0.5 &       0.00018 &          leaky\_relu &         32 &                      256 &                      128 &                      128 &                     512 \\
CM\_0.5\_1.0 &       0.00018 &          leaky\_relu &         32 &                      256 &                      128 &                      128 &                     512 \\
CM\_1.0\_1.0 &       0.00025 &                relu &         64 &                      128 &                      256 &                       64 &                    1024 \\
VRP\_4\_L    &       0.00018 &          leaky\_relu &         32 &                      256 &                      128 &                      128 &                     512 \\
VRP\_10\_L   &       0.00396 &                relu &         32 &                      256 &                       64 &                       64 &                    1024 \\
VRP\_15\_L   &       0.00512 &                relu &        128 &                      256 &                      256 &                       64 &                    1024 \\
VRP\_50\_L   &       0.00326 &          leaky\_relu &         16 &                       64 &                       64 &                      128 &                     512 \\
VRP\_4\_H    &       0.00289 &          leaky\_relu &         64 &                       64 &                      128 &                      128 &                    1024 \\
VRP\_10\_H   &       0.00748 &                relu &        128 &                       64 &                      128 &                       64 &                    1024 \\
VRP\_15\_H   &       0.00526 &                relu &         64 &                      128 &                       64 &                      128 &                     512 \\
VRP\_50\_H   &       0.00512 &                relu &        128 &                      256 &                      256 &                       64 &                    1024 \\
\bottomrule
\end{tabular}}
\caption{Supervised learning best configurations from random search.}
\label{tab:sl_best_config}
\end{table*}

\subsection{DQN Model Selection}
\label{app:MD:DQN}

For DQN, we perform random search over the hyperparameters for both the linear and set-based value function estimators.  The hyperparameters  for which a search is conducted are listed in Table~\ref{tab:random_search_space_dqn}.  We omit the description of the batch size, learning rate, and activation function as they are standard.  One-hot period is a binary decision, which if true encodes the period as a one-hot vector in the state representation, and if false encodes the period as an integer.  Double-DQN is a binary decision which if true uses double Q-learning \citep{van2016deep}, and if false uses standard Q-learning.  Hidden dimensions are only applicable to DQN-L and indicate the hidden dimensions for the value function approximator.  Embed hidden dimension \{1,2,3\} and Concat hidden dimension are the dimensions for the permutation-invariant structure of DQN-S, and are similar to those defined in Appendix~\ref{app:MD:SL}.

\begin{table*}[h!]\centering\resizebox{0.9\textwidth}{!}{
\begin{tabular}{l|cc}
\toprule
{Hyper parameter name}                & {DQN-L}                   & {DQN-S}           \\
\midrule
    Batch size                  & $\{16, 32, 64, 128\}$     & $\{16, 32, 64, 128\}$      \\
    One-hot period              & \{True, False\}           & \{True, False\}  \\
    Learning rate               & $[1e^{-5},1e^{-2}]$       & $[1e^{-5},1e^{-2}]$    \\
    Activation Function         & \{ReLU, LeakyReLU\}       & \{ReLU, LeakyReLU\}   \\
    \# Episodes                 & $15000$                   &  $15000$    \\
    Double DQN                  & \{True, False\}           & \{True, False\} \\
    Hidden dimensions           & $\{[128], [256], [512], [128, 128], [256, 256], [512, 512],  [128, 128, 128],  [256, 256, 256], [512, 512, 512]\}$    & -      \\
    Embed hidden dimension 1    & -                             & $\{64, 128, 256, 512\}$   \\
    Embed hidden dimension 2    & -                             & $\{16, 32, 64, 128\}$  \\
    Embed hidden dimension 3    & -                             & $\{8, 16, 32, 64\}$  \\
    Concat hidden dimension     & -                              & $\{64, 128, 256, 512\}$  \\
\bottomrule
\end{tabular}}
\caption{Random search hyperparameter space for DQN-L and DQN-S.  For values in \{\}, we sample with equal probability for each discrete choice.  For values in [], we sample a uniform distribution within the given bounds.  For single values, we keep it fixed across all configurations.  Note that for the hidden dimension, the value contained in [] are a list of the hidden dimension. A value of - indicates that hyperparameter is not applicable for the given model.}
\label{tab:random_search_space_dqn}
\end{table*}

\begin{figure}[h!]
    \centering
    \includegraphics[width=1.0\textwidth]{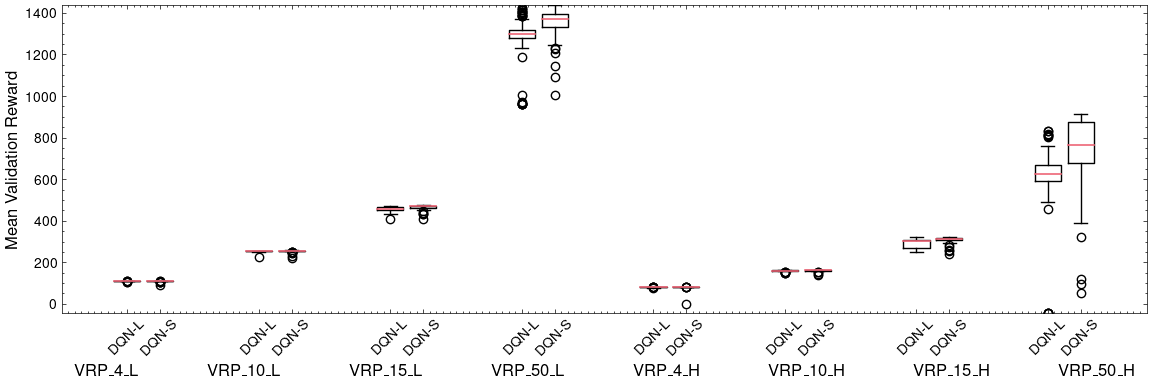}
    \caption[DiL box plots of reward over random search configurations. ]%
    {{\small DiL box plots of reward over random search configurations. }} 
    \label{fig:vrp_dqn_config}
\end{figure}

\begin{figure}[h!]
    \centering
    \includegraphics[width=0.5\textwidth]{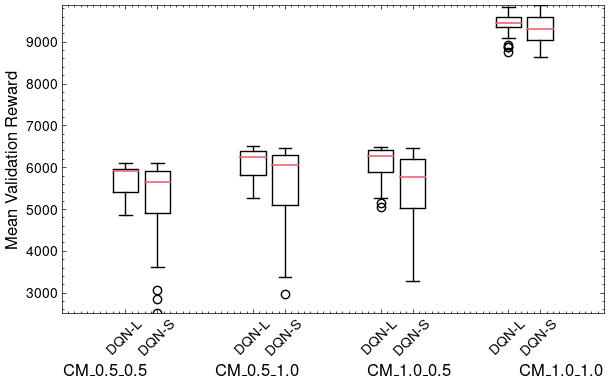}
    \caption[ACM box plots of reward over random search configurations. ]%
    {{\small ACM box plots of reward over random search configurations. }} 
    \label{fig:cm_dqn_config}
\end{figure}

Within the search space specified in Table~\ref{tab:random_search_space_dqn}, we perform random search over 100 configurations for each of DQN-L and DQN-S.  The configuration which achieves the best mean reward over 100 validation trajectories is reported in Tables~\ref{tab:dqn-l_best_config} and \ref{tab:dqn-s_best_config} for the DQN-L and SQN-S, respectively. Box plots of the variations in performance over the various hyperparameter configurations appear in Figures~\ref{fig:vrp_dqn_config}-\ref{fig:cm_dqn_config}.  Here, we report the best mean reward over the same 100 validation instances obtained throughout the training for each configuration.  We conclude that hyperparameter search leads to marginal improvements. For most hyperparameter configurations, high-quality control is achieved.

\begin{table*}[t]\centering\resizebox{1.0\textwidth}{!}{
\begin{tabular}{l|rrrrrrr}
\toprule
{} &  One-hot period & Learning rate & Activation function & Batch size &  Dueling DQN &  Double DQN & Hidden dimensions \\
\midrule
VRP\_4\_L    &           False &       0.00092 &                ReLU &         32 &        False &       False &             [512] \\
VRP\_10\_L   &           False &       0.00224 &           LeakyReLU &         64 &        False &        True &   [256, 256, 256] \\
VRP\_15\_L   &            True &       0.00167 &           LeakyReLU &         32 &        False &       False &        [256, 256] \\
VRP\_50\_L   &           False &       0.00136 &           LeakyReLU &         16 &        False &        True &   [128, 128, 128] \\
VRP\_4\_H    &            True &        0.0013 &           LeakyReLU &         64 &        False &       False &        [256, 256] \\
VRP\_10\_H   &           False &       0.00252 &                ReLU &         16 &        False &       False &   [512, 512, 512] \\
VRP\_15\_H   &            True &       0.00267 &                ReLU &         16 &        False &       False &             [256] \\
VRP\_50\_H   &           False &       0.00508 &           LeakyReLU &         64 &        False &       False &             [512] \\
CM\_0.5\_0.5 &            True &       0.00131 &                ReLU &         32 &        False &       False &        [512, 512] \\
CM\_1.0\_0.5 &            True &       0.00167 &           LeakyReLU &         32 &        False &       False &        [256, 256] \\
CM\_0.5\_1.0 &            True &       0.00118 &           LeakyReLU &         32 &        False &       False &             [256] \\
CM\_1.0\_1.0 &           False &       0.00092 &                ReLU &         32 &        False &       False &             [512] \\
\bottomrule
\end{tabular}}
\caption{DQN-L best configurations from random search.}
\label{tab:dqn-l_best_config}
\end{table*}

\begin{table*}[t]\centering\resizebox{1.0\textwidth}{!}{
\begin{tabular}{l|rrrrrrrrrr}
\toprule
{} &  One-hot period & Learning rate & Activation function & Batch size & Dueling DQN &  Double DQN & Embed hidden dimensions 1 & Embed hidden dimensions 2 & Embed hidden dimensions 3 & Concat hidden dimensions 1 \\
\midrule
VRP\_4\_L    &           False &        0.0021 &           LeakyReLU &         64 &           0 &       False &                        64 &                       512 &                       512 &                        512 \\
VRP\_10\_L   &           False &       0.00258 &                ReLU &         64 &           0 &        True &                       256 &                        64 &                       512 &                        128 \\
VRP\_15\_L   &            True &       0.00036 &           LeakyReLU &         16 &           0 &       False &                       256 &                        64 &                       128 &                       1024 \\
VRP\_50\_L   &            True &       0.00182 &           LeakyReLU &        128 &           0 &       False &                       128 &                       256 &                        64 &                        256 \\
VRP\_4\_H    &            True &       0.00131 &                ReLU &         32 &           0 &       False &                       128 &                        64 &                       128 &                        512 \\
VRP\_10\_H   &            True &        0.0033 &                ReLU &         32 &           0 &       False &                        64 &                       128 &                        64 &                        512 \\
VRP\_15\_H   &            True &        0.0013 &           LeakyReLU &         64 &           0 &       False &                        64 &                       128 &                       128 &                        512 \\
VRP\_50\_H   &            True &        0.0013 &           LeakyReLU &         64 &           0 &       False &                        64 &                       128 &                       128 &                        512 \\
CM\_0.5\_0.5 &           False &       0.00227 &           LeakyReLU &        128 &           0 &       False &                       128 &                       256 &                       128 &                        128 \\
CM\_1.0\_0.5 &           False &       0.00227 &           LeakyReLU &        128 &           0 &       False &                       128 &                       256 &                       128 &                        128 \\
CM\_0.5\_1.0 &           False &       0.00136 &           LeakyReLU &         16 &           0 &        True &                       256 &                        64 &                       512 &                        128 \\
CM\_1.0\_1.0 &           False &       0.00258 &                ReLU &         64 &           0 &        True &                       256 &                        64 &                       512 &                        128 \\
\bottomrule
\end{tabular}}
\caption{DQN-S best configurations from random search.}
\label{tab:dqn-s_best_config}
\end{table*}

\section{Sampling EoHPs in Airline Cargo Management}
\label{app:N_FINAL_OBS}

\begin{figure}[h!]
    \centering
    \includegraphics[width=0.5\textwidth]{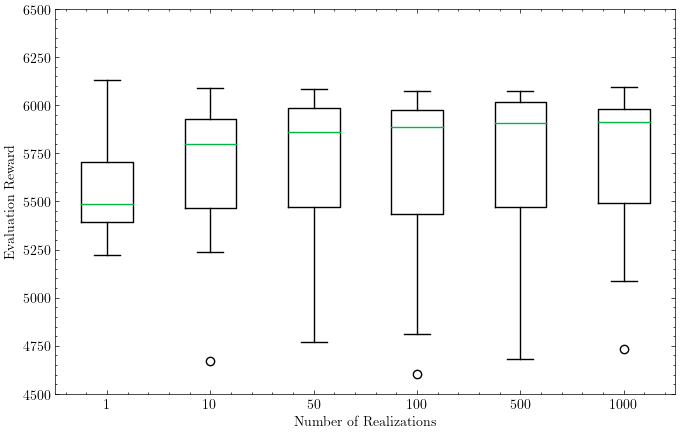}
    \caption[Box plots of evaluation reward for a varying \# of realizations. ]%
    {{\small Box plots of evaluation reward for a varying \# of realizations. }} 
    \label{fig:n_final_obs_var}
\end{figure}

In this section, we report DQN-L results while varying the number of realizations for the EoHP in the ACM problem, specifically for CM\_0.5\_0.5.  Figure~\ref{fig:n_final_obs_var} reports the best evaluation reward obtained throughout DQN-L training of 50 hyperparameter configurations (same configurations for each number of realizations).  From the plot, we can see that the median evaluation reward is notably worse with only a single realization.  Furthermore, as we are only sampling a single realization of the EoHP for each trajectory, we expect that using the evaluation reward may overfit the specific realization rather than being robust over an average of many realizations.  For the experiments in the paper, we use 50 realizations as the distributions are similar for larger numbers (see figure).
Note, however, that our methodology can easily accommodate larger realizations with negligible increases in simulation time.

\section{Reinforcement Learning Training}
\label{app:DQN_CURVES}

\begin{figure*}
    \centering
    \begin{subfigure}[b]{0.75\textwidth}
        \centering
        \includegraphics[width=\textwidth]{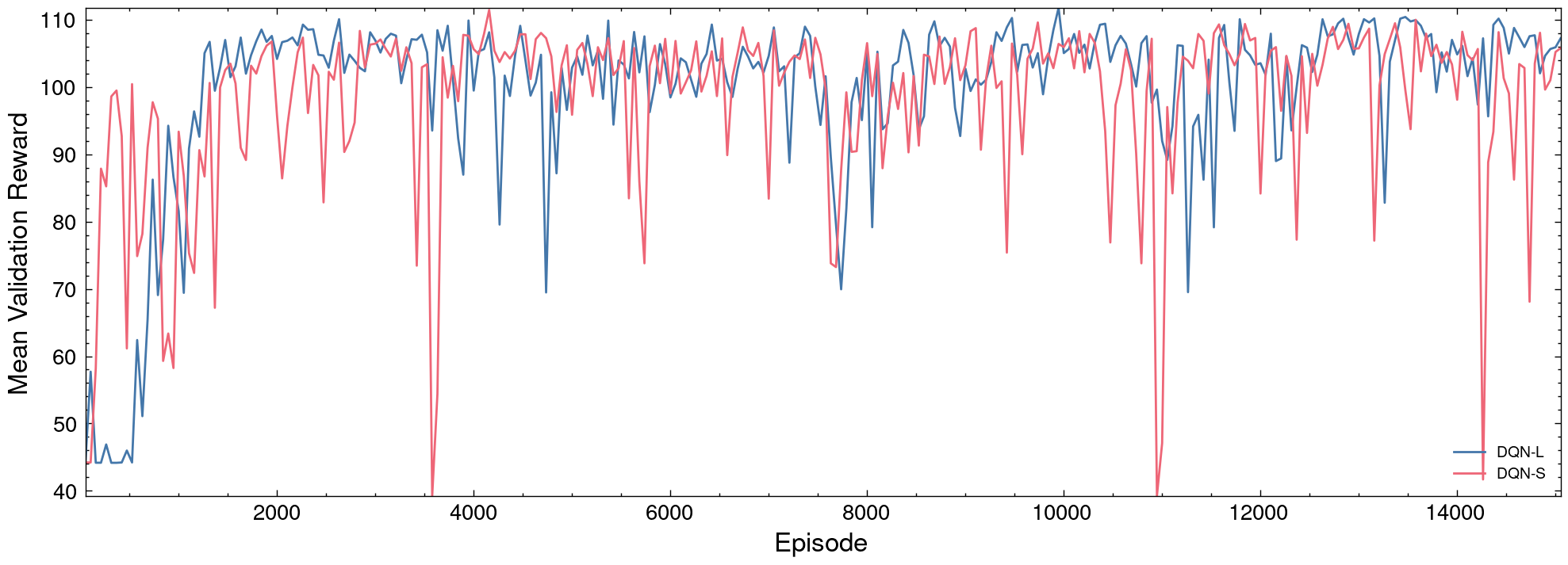}
        \caption[VRP\_4\_L]%
        {{\small VRP\_4\_L}}    
        \label{fig:vrp_4_L_dqn}
    \end{subfigure}
    \hfill
    \begin{subfigure}[b]{0.75\textwidth}  
        \centering 
        \includegraphics[width=\textwidth]{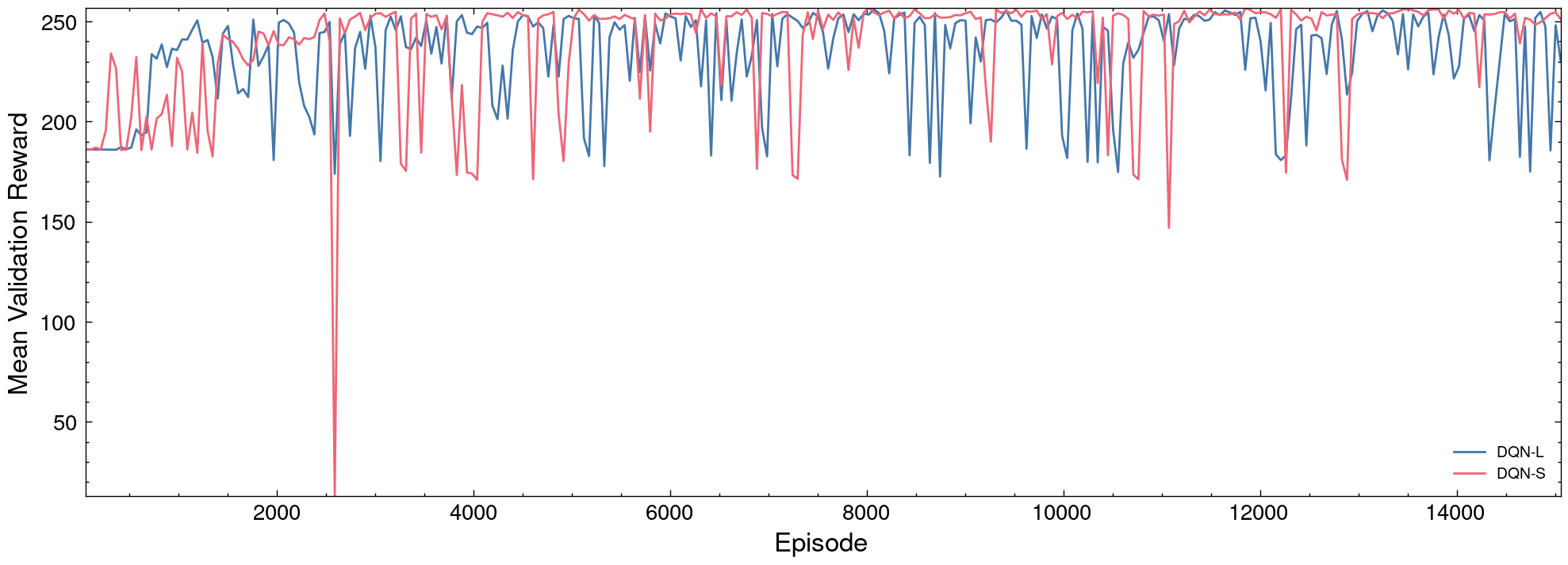}
        \caption[VRP\_10\_L]%
        {{\small VRP\_10\_L}}    
        \label{fig:vrp_10_L_dqn}
    \end{subfigure}
    \vskip\baselineskip
    \begin{subfigure}[b]{0.75\textwidth}   
        \centering 
        \includegraphics[width=\textwidth]{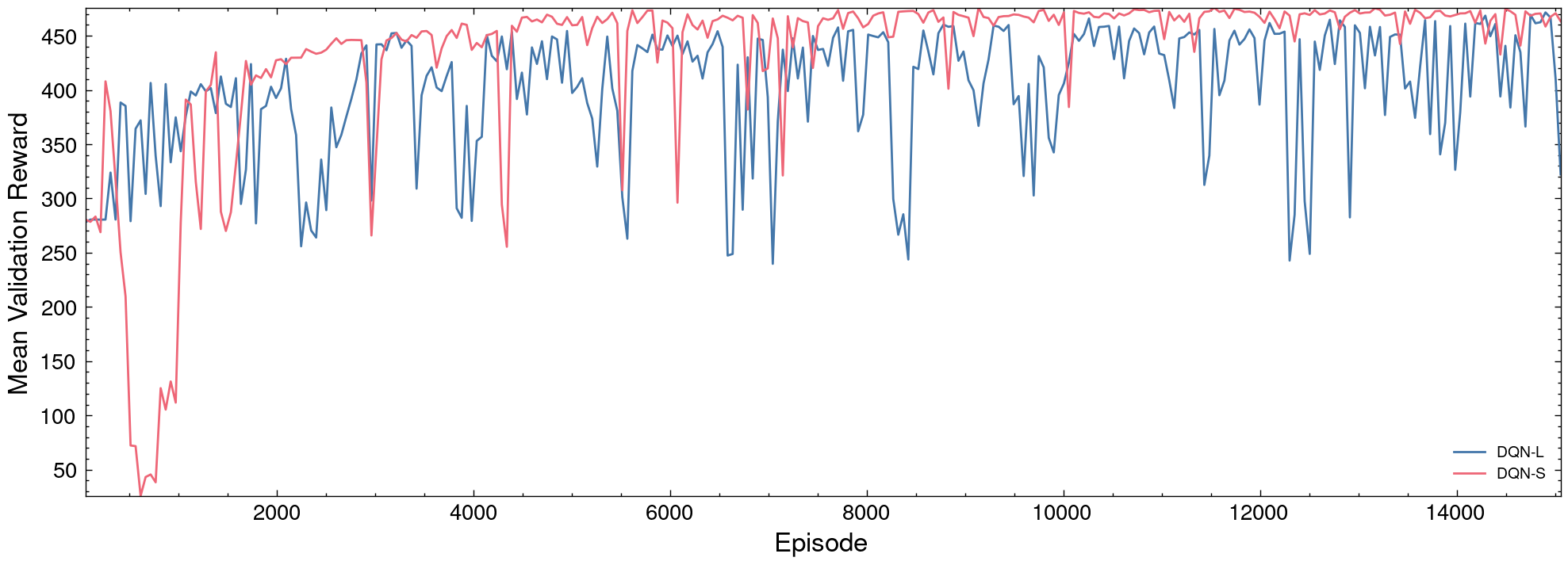}
        \caption[VRP\_15\_L]%
        {{\small VRP\_15\_L}}    
        \label{fig:vrp_15_L_dqn}
    \end{subfigure}
    \hfill
    \begin{subfigure}[b]{0.75\textwidth}   
        \centering 
        \includegraphics[width=\textwidth]{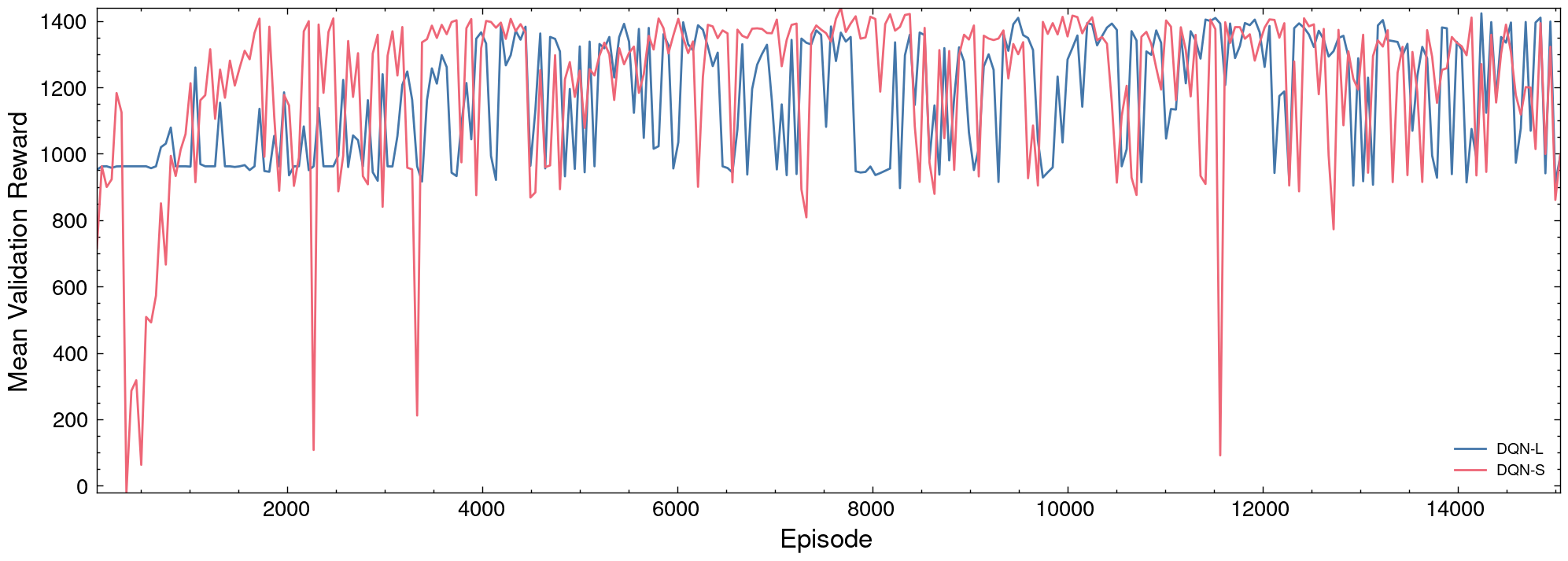}
        \caption[VRP\_50\_L]%
        {{\small VRP\_50\_L}}    
        \label{fig:vrp_50_L_dqn}
    \end{subfigure}
    \caption[DQN training curves for VRP\_X\_L.]
    {{\small DQN training curves for VRP\_X\_L.}} 
    \label{fig:vrp_x_L_dqn}
\end{figure*}

\begin{figure*}
    \centering
    \begin{subfigure}[b]{0.75\textwidth}
        \centering
        \includegraphics[width=\textwidth]{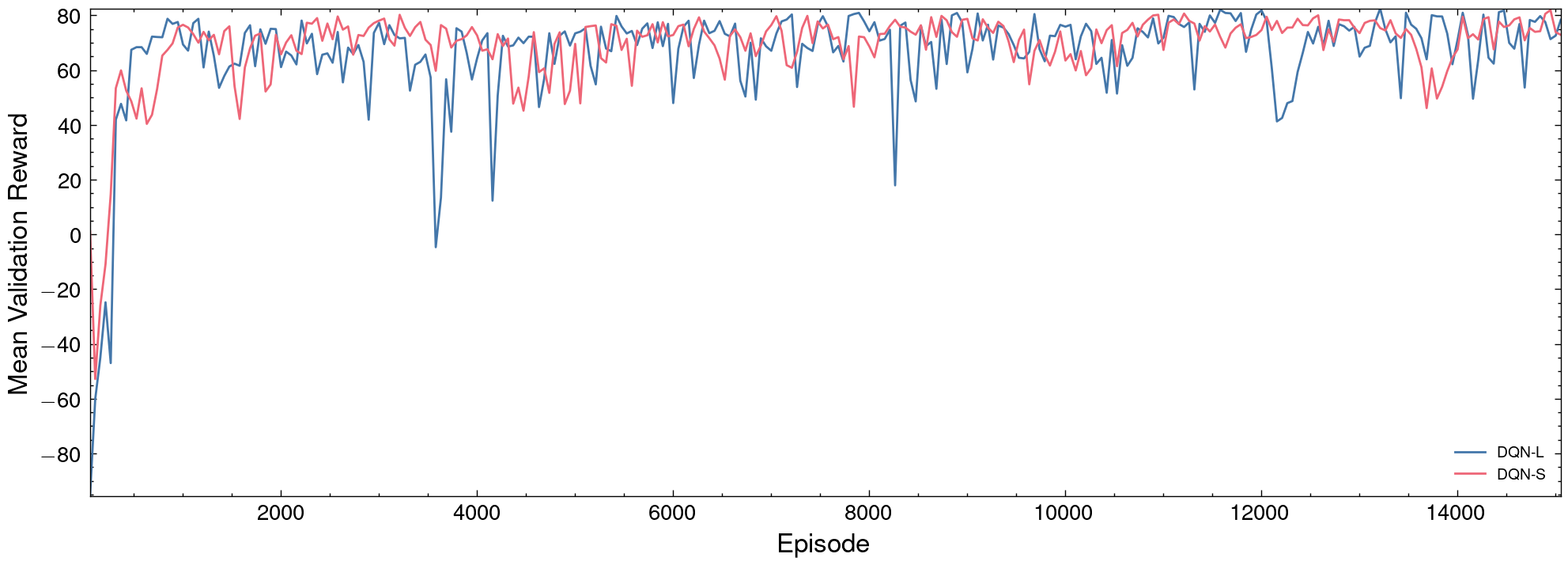}
        \caption[VRP\_4\_H]%
        {{\small VRP\_4\_H}}    
        \label{fig:vrp_4_H_dqn}
    \end{subfigure}
    \hfill
    \begin{subfigure}[b]{0.75\textwidth}  
        \centering 
        \includegraphics[width=\textwidth]{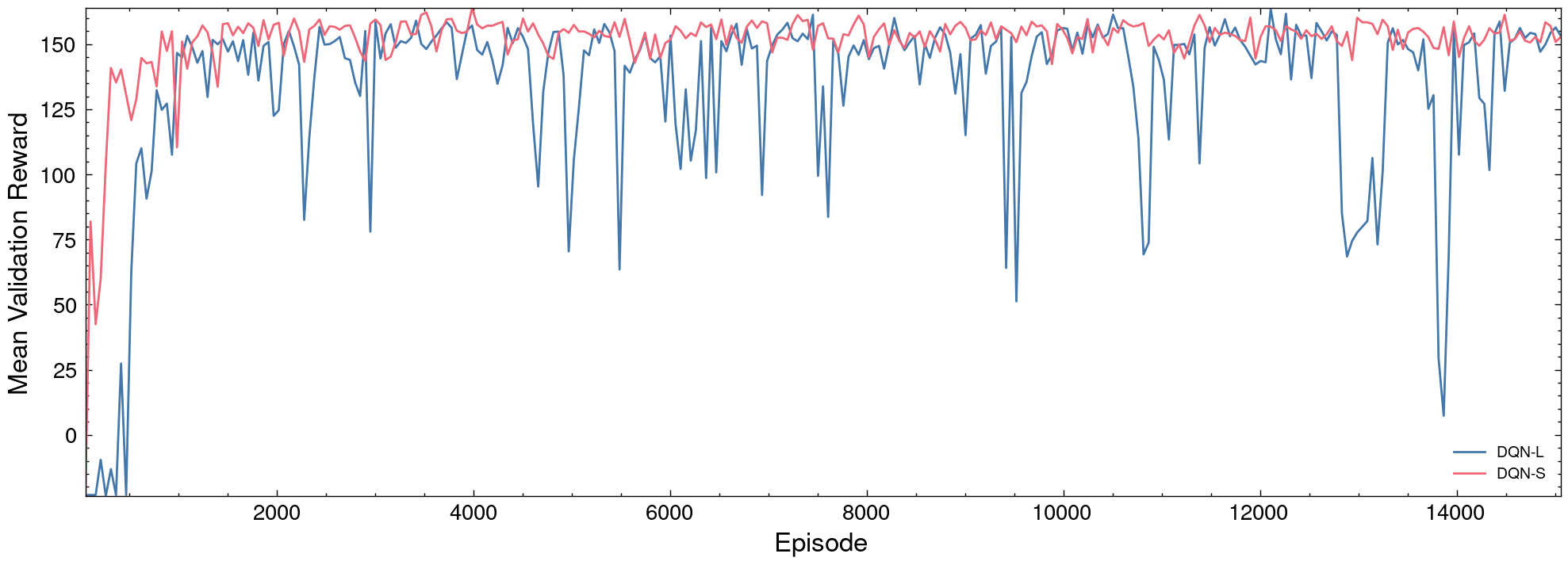}
        \caption[VRP\_10\_H]%
        {{\small VRP\_10\_H}}    
        \label{fig:vrp_10_H_dqn}
    \end{subfigure}
    \vskip\baselineskip
    \begin{subfigure}[b]{0.75\textwidth}   
        \centering 
        \includegraphics[width=\textwidth]{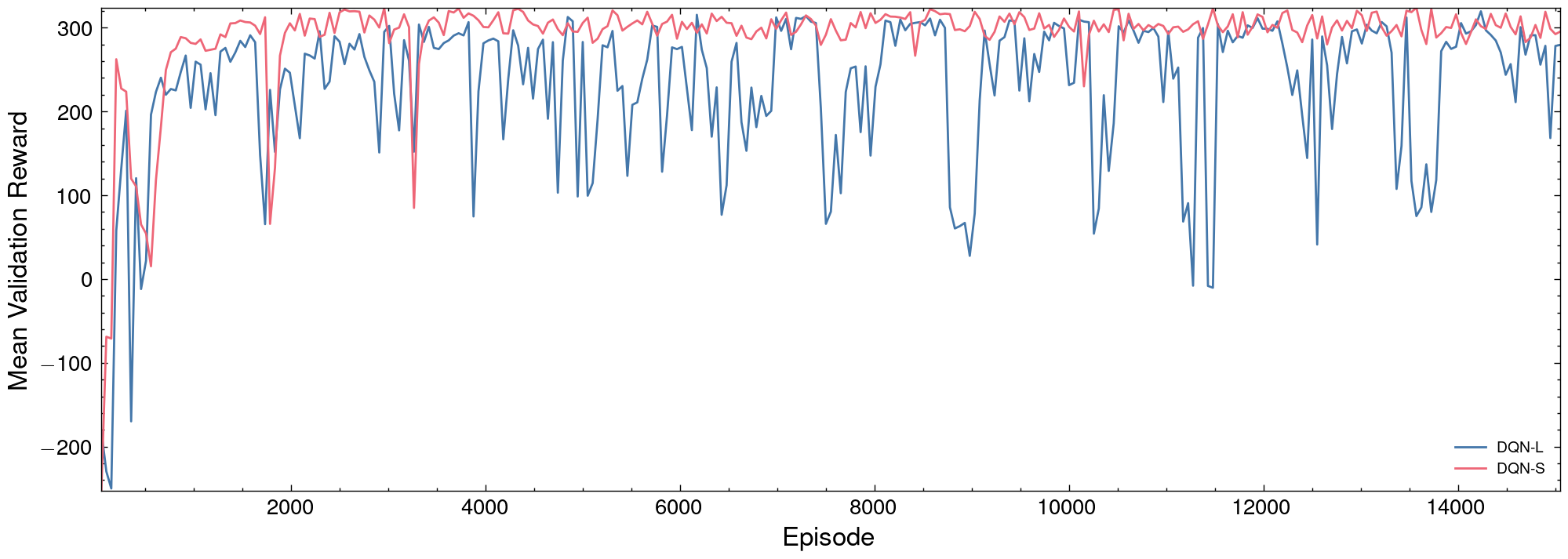}
        \caption[VRP\_15\_H]%
        {{\small VRP\_15\_H}}    
        \label{fig:vrp_15_H_dqn}
    \end{subfigure}
    \hfill
    \begin{subfigure}[b]{0.75\textwidth}   
        \centering 
        \includegraphics[width=\textwidth]{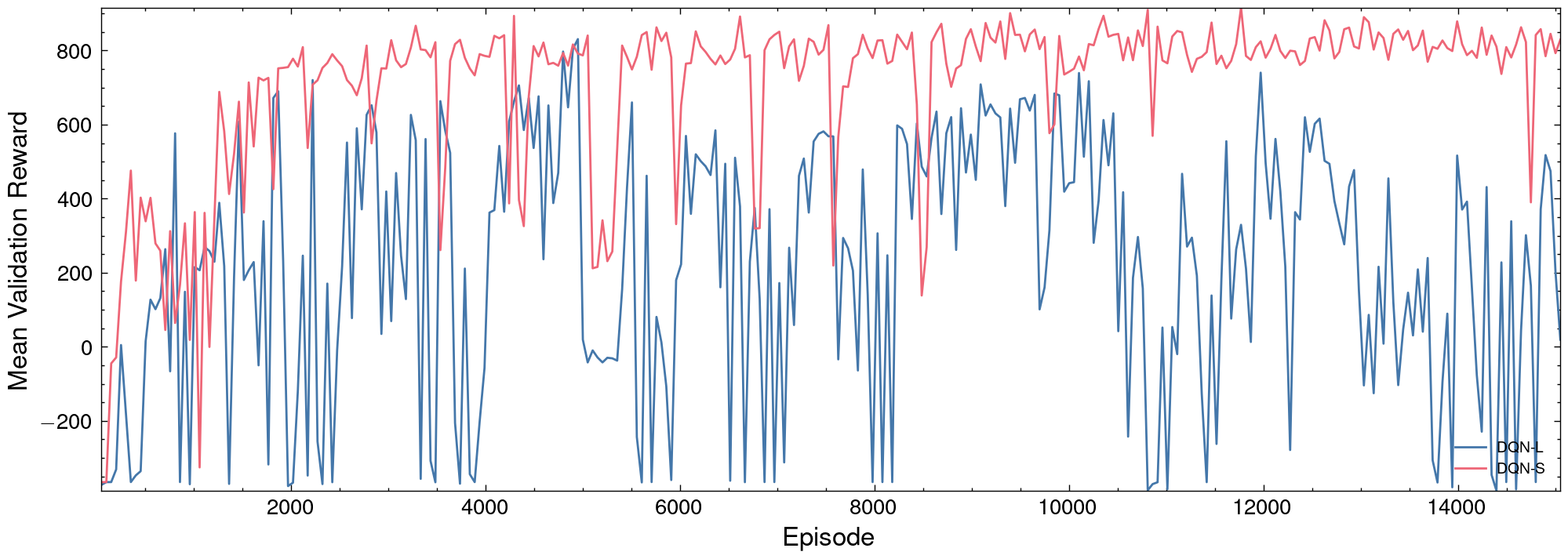}
        \caption[VRP\_50\_H]%
        {{\small VRP\_50\_H}}    
        \label{fig:vrp_50_H_dqn}
    \end{subfigure}
    \caption[ DQN training curves for VRP\_X\_H. ]
    {{\small DQN training curves for VRP\_X\_H.}} 
    \label{fig:vrp_x_H_dqn}
\end{figure*}

\begin{figure*}
    \centering
    \begin{subfigure}[b]{0.75\textwidth}
        \centering
        \includegraphics[width=\textwidth]{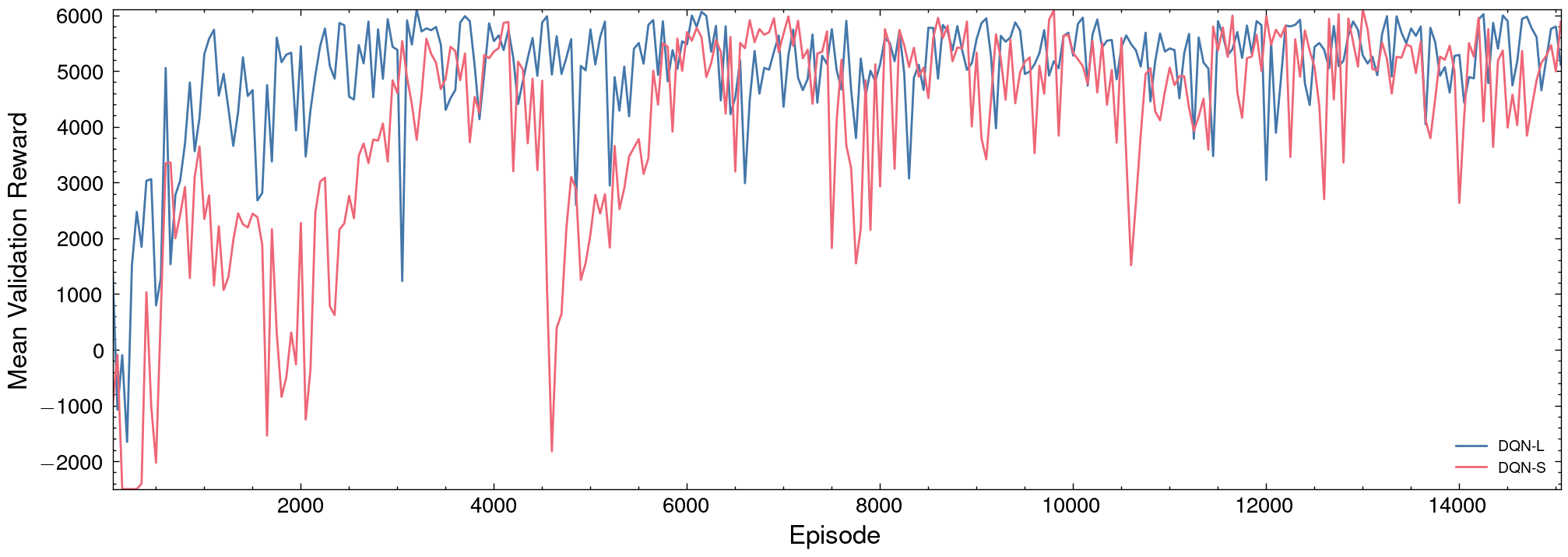}
        \caption[CM\_0.5\_0.5]%
        {{\small CM\_0.5\_0.5}}    
        \label{fig:cm_2_2_dqn}
    \end{subfigure}
    \hfill
    \begin{subfigure}[b]{0.75\textwidth}  
        \centering 
        \includegraphics[width=\textwidth]{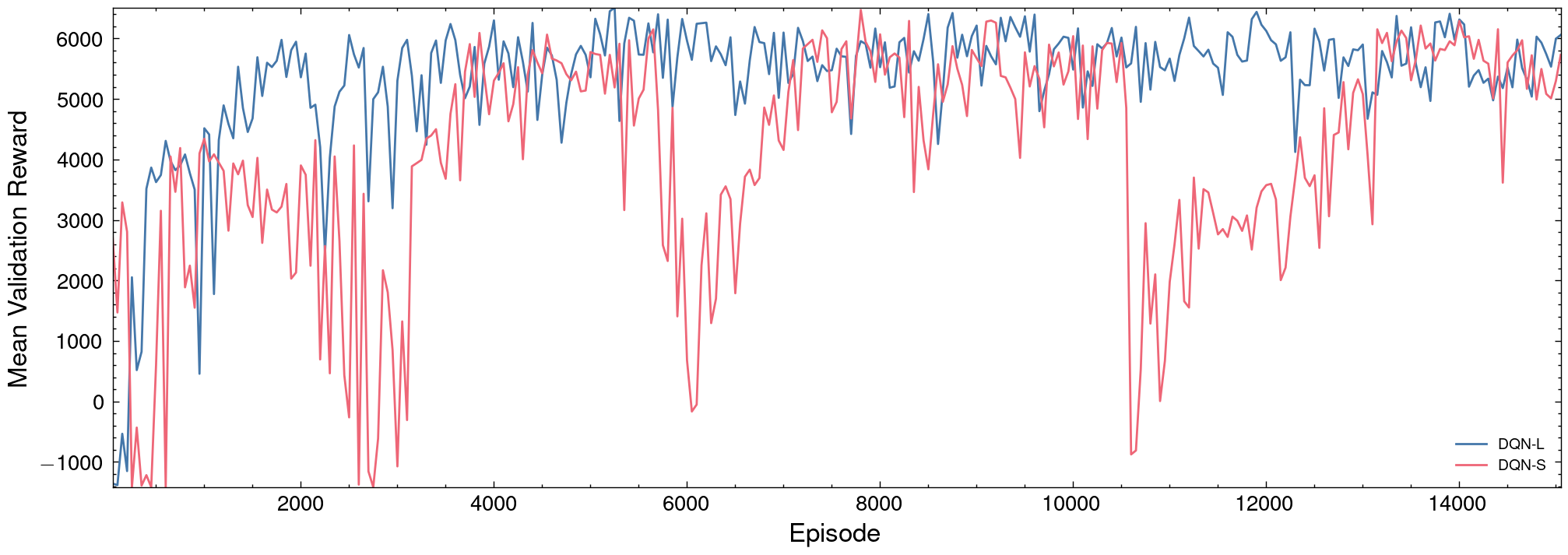}
        \caption[CM\_0.5\_1.0]%
        {{\small CM\_0.5\_1.0}}    
        \label{fig:cm_2_1_dqn}
    \end{subfigure}
    \vskip\baselineskip
    \begin{subfigure}[b]{0.75\textwidth}   
        \centering 
        \includegraphics[width=\textwidth]{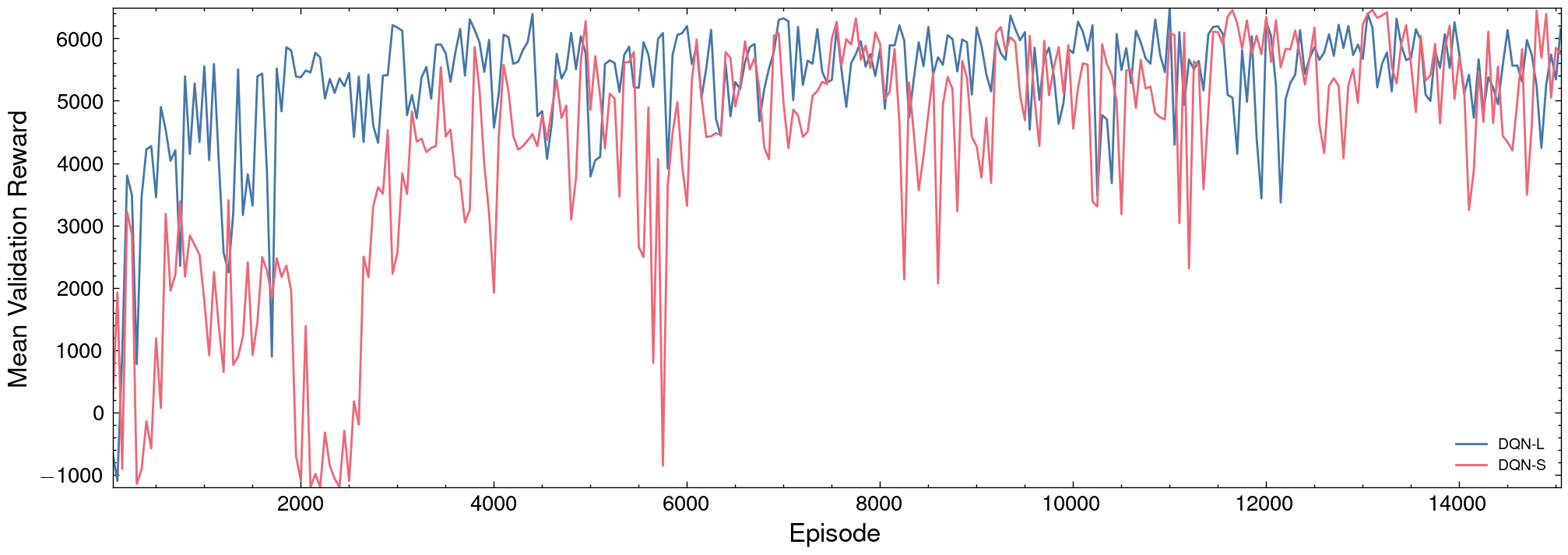}
        \caption[CM\_1.0\_0.5]%
        {{\small CM\_1.0\_0.5}}    
        \label{fig:cm_1_2_dqn}
    \end{subfigure}
    \hfill
    \begin{subfigure}[b]{0.75\textwidth}   
        \centering 
        \includegraphics[width=\textwidth]{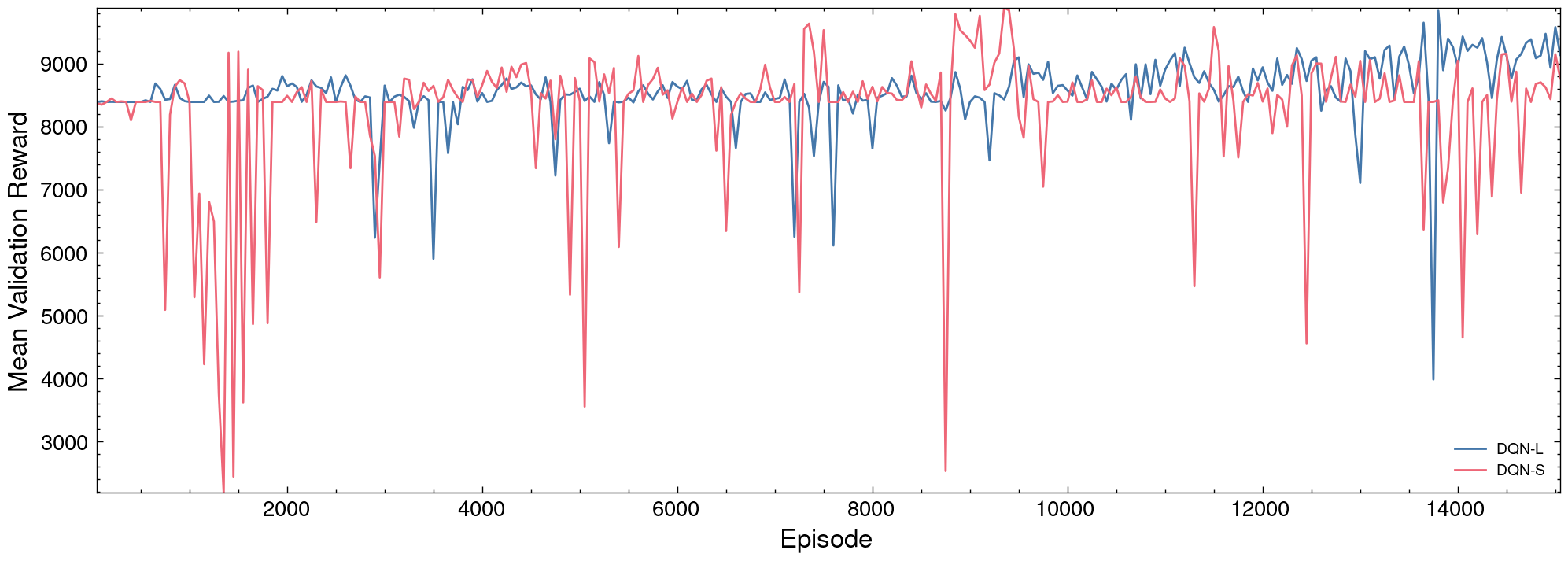}
        \caption[CM\_1.0\_1.0]%
        {{\small CM\_1.0\_1.0}}    
        \label{fig:cm_1_1_dqn}
    \end{subfigure}
    \caption[ DQN training curves for CM\_X\_Y. ]
    {{\small DQN training curves for CM\_X\_Y.}} 
    \label{fig:cm_dqn}
\end{figure*}

For completeness, we provide in Figures~\ref{fig:vrp_x_L_dqn} to \ref{fig:cm_dqn} the training curves for the best DQN configuration for each instance.  These curves report the mean reward over the 100 validation trajectories for each of the instances over the 15,000 episodes.  We notice a relatively high variance between episodes.  This is expected, given that we are only plotting the reward curve for 1 seed, rather than multiple which is standard in RL.  The curves indicate that convergence occurs in most cases and that the difference between DQN-L and DQN-S is often negligible, with the exception of the larger VRP instances.  It is useful to note that convergence occurs typically in less than 5,000 episodes so that training with 15,000 may not be required.

\clearpage

\bibliographystyle{plainnat}  
\bibliography{references}

\begin{thebibliography}{53}
\providecommand{\natexlab}[1]{#1}
\providecommand{\url}[1]{\texttt{#1}}
\expandafter\ifx\csname urlstyle\endcsname\relax
  \providecommand{\doi}[1]{doi: #1}\else
  \providecommand{\doi}{doi: \begingroup \urlstyle{rm}\Url}\fi

\bibitem[Accorsi and Vigo(2021)]{accorsi2021fast}
Luca Accorsi and Daniele Vigo.
\newblock A fast and scalable heuristic for the solution of large-scale
  capacitated vehicle routing problems.
\newblock \emph{Transportation Science}, 55\penalty0 (4):\penalty0 832--856,
  2021.

\bibitem[Adelman(2007)]{adelman2007dynamic}
Daniel Adelman.
\newblock Dynamic bid prices in revenue management.
\newblock \emph{Operations Research}, 55\penalty0 (4):\penalty0 647--661, 2007.

\bibitem[Amaruchkul et~al.(2007)Amaruchkul, Cooper, and
  Gupta]{amaruchkul2007single}
Kannapha Amaruchkul, William~L Cooper, and Diwakar Gupta.
\newblock Single-leg air-cargo revenue management.
\newblock \emph{Transportation Science}, 41\penalty0 (4):\penalty0 457--469,
  2007.

\bibitem[Barz and Gartner(2016)]{barz2016air}
Christiane Barz and Daniel Gartner.
\newblock Air cargo network revenue management.
\newblock \emph{Transportation Science}, 50\penalty0 (4):\penalty0 1206--1222,
  2016.

\bibitem[Bengio et~al.(2021)Bengio, Lodi, and Prouvost]{bengio2020machine}
Yoshua Bengio, Andrea Lodi, and Antoine Prouvost.
\newblock Machine learning for combinatorial optimization: {A} methodological
  tour d’horizon.
\newblock \emph{European Journal of Operational Research}, 290\penalty0
  (2):\penalty0 405--421, 2021.

\bibitem[Bertsekas(2019)]{bertsekas2019reinforcement}
Dimitri~P Bertsekas.
\newblock \emph{Reinforcement learning and optimal control}.
\newblock Athena Scientific Belmont, MA, 2019.

\bibitem[Bertsimas and Popescu(2003)]{bertsimas2003revenue}
Dimitris Bertsimas and Ioana Popescu.
\newblock Revenue management in a dynamic network environment.
\newblock \emph{Transportation Science}, 37\penalty0 (3):\penalty0 257--277,
  2003.

\bibitem[Bondoux et~al.(2020)Bondoux, Nguyen, Fiig, and
  Acuna-Agost]{bondoux2020reinforcement}
Nicolas Bondoux, Anh~Quan Nguyen, Thomas Fiig, and Rodrigo Acuna-Agost.
\newblock Reinforcement learning applied to airline revenue management.
\newblock \emph{Journal of Revenue and Pricing Management}, 19\penalty0
  (5):\penalty0 332--348, 2020.

\bibitem[Brandt and Nickel(2019)]{nickel}
Felix Brandt and Stefan Nickel.
\newblock The air cargo load planning problem -- {A} consolidated problem
  definition and literature review on related problems.
\newblock \emph{European Journal of Operational Research}, 275:\penalty0
  399--410, 2019.

\bibitem[Breiman(2001)]{breiman2001random}
Leo Breiman.
\newblock Random forests.
\newblock \emph{Machine Learning}, 45\penalty0 (1):\penalty0 5--32, 2001.

\bibitem[Brown et~al.(2010)Brown, Smith, and Sun]{brown2010information}
David~B Brown, James~E Smith, and Peng Sun.
\newblock Information relaxations and duality in stochastic dynamic programs.
\newblock \emph{Operations Research}, 58\penalty0 (4-part-1):\penalty0
  785--801, 2010.

\bibitem[Brown et~al.(2022)Brown, Smith, et~al.]{brown2022information}
David~B Brown, James~E Smith, et~al.
\newblock Information relaxations and duality in stochastic dynamic programs: A
  review and tutorial.
\newblock \emph{Foundations and Trends in Optimization}, 5\penalty0
  (3):\penalty0 246--339, 2022.

\bibitem[Buconiu et~al.(2018)Buconiu, {de Bruin}, Tolic, Kober, and
  Palunko]{BusonEtAl18}
Lucian Buconiu, Tim {de Bruin}, Domagoj Tolic, Jens Kober, and Ivana Palunko.
\newblock Reinforcement learning for control: Performance, stability, and deep
  approximators.
\newblock \emph{Annual Reviews in Control}, 46:\penalty0 8--28, 2018.

\bibitem[Caprara and Toth(2001)]{caprara2001lower}
Alberto Caprara and Paolo Toth.
\newblock Lower bounds and algorithms for the 2-dimensional vector packing
  problem.
\newblock \emph{Discrete Applied Mathematics}, 111\penalty0 (3):\penalty0
  231--262, 2001.

\bibitem[Dai et~al.(2022)Dai, Xue, Syed, Schuurmans, and Dai]{dai2022neural}
Hanjun Dai, Yuan Xue, Zia Syed, Dale Schuurmans, and Bo~Dai.
\newblock Neural stochastic dual dynamic programming.
\newblock In \emph{International Conference on Learning Representations}, 2022.

\bibitem[Donti et~al.(2021)Donti, Roderick, Fazlyab, and
  Kolter]{donti2021enforcing}
Priya~L. Donti, Melrose Roderick, Mahyar Fazlyab, and J~Zico Kolter.
\newblock Enforcing robust control guarantees within neural network policies.
\newblock In \emph{In International Conference on Learning Representations},
  2021.

\bibitem[Dumouchelle et~al.(2022)Dumouchelle, Patel, Khalil, and
  Bodur]{dumouchelle2022neur2sp}
Justin Dumouchelle, Rahul Patel, Elias~B Khalil, and Merve Bodur.
\newblock Neur2{SP}: Neural two-stage stochastic programming.
\newblock \emph{Advances in Neural Information Processing Systems}, 35, 2022.

\bibitem[Finn et~al.(2017)Finn, Yu, Fu, Abbeel, and
  Levine]{finn2017generalizing}
Chelsea Finn, Tianhe Yu, Justin Fu, Pieter Abbeel, and Sergey Levine.
\newblock Generalizing skills with semi-supervised reinforcement learning.
\newblock In \emph{International Conference on Learning Representations}, 2017.

\bibitem[Fischetti and Fraccaro(2017)]{fischetti2017using}
Martina Fischetti and Marco Fraccaro.
\newblock Using {OR} + {AI} to predict the optimal production of offshore wind
  parks: a preliminary study.
\newblock In \emph{International Conference on Optimization and Decision
  Science}, pages 203--211. Springer, 2017.

\bibitem[Giallombardo et~al.(2022)Giallombardo, Guerriero, and Miglionico]{ggm}
Giovanni Giallombardo, Francesca Guerriero, and Giovanna Miglionico.
\newblock Profit maximization via capacity control for distribution logistics
  problems.
\newblock \emph{Computers \& Industrial Engineering}, 171:\penalty0 108466,
  2022.

\bibitem[Gosavi(2009)]{Gosavi09}
Abhijit Gosavi.
\newblock Reinforcement learning: A tutorial survey and recent advances.
\newblock \emph{INFORMS Journal on Computing}, 21\penalty0 (2):\penalty0
  178--192, 2009.

\bibitem[Gosavi et~al.(2002)Gosavi, Bandla, and Das]{gosavii2002reinforcement}
Abhuit Gosavi, Naveen Bandla, and Tapas~K. Das.
\newblock A reinforcement learning approach to a single leg airline revenue
  management problem with multiple fare classes and overbooking.
\newblock \emph{{IIE} Transactions}, 34\penalty0 (9):\penalty0 729--742, 2002.

\bibitem[{Gurobi Optimization, LLC}(2021)]{gurobi}
{Gurobi Optimization, LLC}.
\newblock {Gurobi Optimizer Reference Manual}, 2021.
\newblock URL \url{https://www.gurobi.com/documentation/9.5/refman/index.html}.

\bibitem[Henderson et~al.(2005)Henderson, Lemon, and
  Georgila]{henderson2005hybrid}
James Henderson, Oliver Lemon, and Kallirroi Georgila.
\newblock Hybrid reinforcement/supervised learning for dialogue policies from
  communicator data.
\newblock In \emph{IJCAI workshop on knowledge and reasoning in practical
  dialogue systems}, 2005.

\bibitem[Kasilingam(1997)]{kasilingam1997air}
Raja~G Kasilingam.
\newblock Air cargo revenue management: Characteristics and complexities.
\newblock \emph{European Journal of Operational Research}, 96\penalty0
  (1):\penalty0 36--44, 1997.

\bibitem[Kastius and Schlosser(2022)]{KastSchl22}
Alexander Kastius and Rainer Schlosser.
\newblock Dynamic pricing under competition using reinforcement learning.
\newblock \emph{Journal of Revenue and Pricing Management}, 21\penalty0
  (1):\penalty0 50--63, 2022.

\bibitem[Kephart and Tesauro(2000)]{KephTesa00}
Jeffrey~O. Kephart and Gerald~J. Tesauro.
\newblock Pseudo-convergent {Q}-learning by competitive pricebots.
\newblock In \emph{Proceedings 17th International Conference Machine Learning},
  pages 463--470. Morgan Kaufmann, 2000.

\bibitem[Kotary et~al.(2021)Kotary, Fioretto, Van~Hentenryck, and
  Wilder]{KotaryEtAl21}
James Kotary, Ferdinando Fioretto, Pascal Van~Hentenryck, and Bryan Wilder.
\newblock End-to-end constrained optimization learning: A survey.
\newblock In \emph{Proceedings of the Thirtieth International Joint Conference
  on Artificial Intelligence, {IJCAI-21}}, pages 4475--4482, 8 2021.

\bibitem[Larsen et~al.(2022{\natexlab{a}})Larsen, Frejinger, Gendron, and
  Lodi]{LarsFrejGendLodi22}
Eric Larsen, Emma Frejinger, Bernard Gendron, and Andrea Lodi.
\newblock Fast continuous and integer {L}-shaped heuristics through supervised
  learning, 2022{\natexlab{a}}.
\newblock arXiv preprint arXiv:2205.00897v2.

\bibitem[Larsen et~al.(2022{\natexlab{b}})Larsen, Lachapelle, Bengio,
  Frejinger, Lacoste-Julien, and Lodi]{LarsEtAl22}
Eric Larsen, S\'{e}bastien Lachapelle, Yoshua Bengio, Emma Frejinger, Simon
  Lacoste-Julien, and Andrea Lodi.
\newblock Predicting tactical solutions to operational planning problems under
  imperfect information.
\newblock \emph{INFORMS Journal on Computing}, 34\penalty0 (1):\penalty0
  227--242, 2022{\natexlab{b}}.

\bibitem[Lawhead and Gosavi(2019)]{lawhead2019bounded}
Ryan~J. Lawhead and Abhijit Gosavi.
\newblock A bounded actor–critic reinforcement learning algorithm applied to
  airline revenue management.
\newblock \emph{Engineering Applications of Artificial Intelligence},
  82:\penalty0 252 -- 262, 2019.

\bibitem[Levin et~al.(2012)Levin, Nediak, and Topaloglu]{levin2012cargo}
Yuri Levin, Mikhail Nediak, and Huseyin Topaloglu.
\newblock Cargo capacity management with allotments and spot market demand.
\newblock \emph{Operations Research}, 60\penalty0 (2):\penalty0 351--365, 2012.

\bibitem[Levina et~al.(2011)Levina, Levin, McGill, and
  Nediak]{levina2011network}
Tatsiana Levina, Yuri Levin, Jeff McGill, and Mikhail Nediak.
\newblock Network cargo capacity management.
\newblock \emph{Operations Research}, 59\penalty0 (4):\penalty0 1008--1023,
  2011.

\bibitem[Liang et~al.(2018)Liang, Liaw, Nishihara, Moritz, Fox, Goldberg,
  Gonzalez, Jordan, and Stoica]{liang2018rllib}
Eric Liang, Richard Liaw, Robert Nishihara, Philipp Moritz, Roy Fox, Ken
  Goldberg, Joseph Gonzalez, Michael Jordan, and Ion Stoica.
\newblock Rllib: Abstractions for distributed reinforcement learning.
\newblock In \emph{International Conference on Machine Learning}, pages
  3053--3062, 2018.

\bibitem[Martello and Toth(1990)]{martello1990knapsack}
Silvano Martello and Paolo Toth.
\newblock \emph{Knapsack problems: {A}lgorithms and computer implementations}.
\newblock John Wiley \& Sons Ltd., USA, 1990.

\bibitem[Mnih et~al.(2015)Mnih, Kavukcuoglu, Silver, Rusu, Veness, Bellemare,
  Graves, Riedmiller, Fidjeland, Ostrovski, Petersen, Beattie, Sadik,
  Antonoglou, King, Kumaran, Wierstra, Legg, and Hassabis]{mnih2015human}
Volodymyr Mnih, Koray Kavukcuoglu, David Silver, Andrei~A. Rusu, Joel Veness,
  Marc~G. Bellemare, Alex Graves, Martin Riedmiller, Andreas~K. Fidjeland,
  Georg Ostrovski, Stig Petersen, Charles Beattie, Amir Sadik, Ioannis
  Antonoglou, Helen King, Dharshan Kumaran, Daan Wierstra, Shane Legg, and
  Demis Hassabis.
\newblock Human-level control through deep reinforcement learning.
\newblock \emph{Nature}, 518\penalty0 (7540):\penalty0 529--533, 2015.

\bibitem[Paszke et~al.(2019)Paszke, Gross, Massa, Lerer, Bradbury, Chanan,
  Killeen, Lin, Gimelshein, Antiga, Desmaison, Kopf, Yang, DeVito, Raison,
  Tejani, Chilamkurthy, Steiner, Fang, Bai, and Chintala]{NEURIPS2019_9015}
Adam Paszke, Sam Gross, Francisco Massa, Adam Lerer, James Bradbury, Gregory
  Chanan, Trevor Killeen, Zeming Lin, Natalia Gimelshein, Luca Antiga, Alban
  Desmaison, Andreas Kopf, Edward Yang, Zachary DeVito, Martin Raison, Alykhan
  Tejani, Sasank Chilamkurthy, Benoit Steiner, Lu~Fang, Junjie Bai, and Soumith
  Chintala.
\newblock Pytorch: {A}n imperative style, high-performance deep learning
  library.
\newblock In \emph{Advances in Neural Information Processing Systems},
  volume~32, pages 8026--8037, 2019.

\bibitem[Pedregosa et~al.(2011)Pedregosa, Varoquaux, Gramfort, Michel, Thirion,
  Grisel, Blondel, Prettenhofer, Weiss, Dubourg, Vanderplas, Passos,
  Cournapeau, Brucher, Perrot, and Duchesnay]{scikit-learn}
F.~Pedregosa, G.~Varoquaux, A.~Gramfort, V.~Michel, B.~Thirion, O.~Grisel,
  M.~Blondel, P.~Prettenhofer, R.~Weiss, V.~Dubourg, J.~Vanderplas, A.~Passos,
  D.~Cournapeau, M.~Brucher, M.~Perrot, and E.~Duchesnay.
\newblock Scikit-learn: {M}achine learning in {P}ython.
\newblock \emph{Journal of Machine Learning Research}, 12:\penalty0 2825--2830,
  2011.

\bibitem[Prates et~al.(2019)Prates, Avelar, Lemos, Lamb, and
  Vardi]{prates2019learning}
Marcelo Prates, Pedro~HC Avelar, Henrique Lemos, Luis~C Lamb, and Moshe~Y
  Vardi.
\newblock Learning to solve {NP}-complete problems: {A} graph neural network
  for decision tsp.
\newblock In \emph{Proceedings of the AAAI Conference on Artificial
  Intelligence}, volume~33, pages 4731--4738, 2019.

\bibitem[Rana and Oliveira(2014)]{RanaEtAl14}
Rupal Rana and Fernando~S. Oliveira.
\newblock Real-time dynamic pricing in a non-stationary environment using
  model-free reinforcement learning.
\newblock \emph{Omega}, 47:\penalty0 116--126, 2014.

\bibitem[Ryu et~al.(2020)Ryu, Chow, Anderson, Tjandraatmadja, and
  Boutilier]{Ryu2020CAQL:}
Moonkyung Ryu, Yinlam Chow, Ross Anderson, Christian Tjandraatmadja, and Craig
  Boutilier.
\newblock {CAQL}: {C}ontinuous action {Q}-learning.
\newblock In \emph{International Conference on Learning Representations}, 2020.

\bibitem[Seo et~al.(2021)Seo, Chang, Cheong, and Baek]{SeoEtAl21}
Dong-Wook Seo, Kyuchang Chang, Taesu Cheong, and Jun-Geol Baek.
\newblock A reinforcement learning approach to distribution-free capacity
  allocation for sea cargo revenue management.
\newblock \emph{Information Sciences}, 571:\penalty0 623--648, 2021.

\bibitem[Shihab et~al.(2019)Shihab, Logemann, Thomas, and
  Wei]{shihab2019autonomous}
Syed Arbab~Mohd Shihab, Caleb Logemann, Deepak-George Thomas, and Peng Wei.
\newblock Autonomous airline revenue management: A deep reinforcement learning
  approach to seat inventory control and overbooking, 2019.
\newblock arXiv preprint arXiv:1902.06824.

\bibitem[Shukla et~al.(2020)Shukla, Kolbeinsson, Marla, and
  Yellepeddi]{ShuklaEtAl20}
Naman Shukla, Arinbjorn Kolbeinsson, Lavanya Marla, and Kartik Yellepeddi.
\newblock From average customer to individual traveler: {A} field experiment in
  airline ancillary pricing.
\newblock \emph{Available at SSRN}, 2020.
\newblock \doi{10.2139/ssrn.3518854}.

\bibitem[Silver et~al.(2016)Silver, Huang, Maddison, Guez, Sifre, van~den
  Driessche, Schrittwieser, Antonoglou, Panneershelvam, Lanctot, Dieleman,
  Grewe, Nham, Kalchbrenner, Sutskever, Lillicrap, Leach, Kavukcuoglu, Graepel,
  and Hassabis]{silver2016mastering}
David Silver, Aja Huang, Chris~J. Maddison, Arthur Guez, Laurent Sifre, George
  van~den Driessche, Julian Schrittwieser, Ioannis Antonoglou, Veda
  Panneershelvam, Marc Lanctot, Sander Dieleman, Dominik Grewe, John Nham, Nal
  Kalchbrenner, Ilya Sutskever, Timothy Lillicrap, Madeleine Leach, Koray
  Kavukcuoglu, Thore Graepel, and Demis Hassabis.
\newblock Mastering the game of {G}o with deep neural networks and tree search.
\newblock \emph{Nature}, 529\penalty0 (7587):\penalty0 484--489, 2016.

\bibitem[Sutton and Barto(2018)]{sutton2018reinforcement}
Richard~S. Sutton and Andrew~G. Barto.
\newblock \emph{Reinforcement Learning: An Introduction}.
\newblock A Bradford Book, Cambridge, MA, USA, 2018.

\bibitem[Talluri and Ryzin(2004)]{talluri2006theory}
Kalyan~T. Talluri and Garrett J.~Van Ryzin.
\newblock \emph{The Theory and Practice of Revenue Management}.
\newblock Springer {US}, 2004.

\bibitem[Toth and Vigo(2002)]{toth2002vehicle}
Paolo Toth and Daniele Vigo.
\newblock \emph{The vehicle routing problem}.
\newblock SIAM, 2002.

\bibitem[Van~Hasselt et~al.(2016)Van~Hasselt, Guez, and Silver]{van2016deep}
Hado Van~Hasselt, Arthur Guez, and David Silver.
\newblock Deep reinforcement learning with double {Q}-learning.
\newblock 30\penalty0 (1), 2016.

\bibitem[Watkins and Dayan(1992)]{watkins1992q}
Christopher~JCH Watkins and Peter Dayan.
\newblock Q-learning.
\newblock \emph{Machine learning}, 8:\penalty0 279--292, 1992.

\bibitem[Williamson(1992)]{williamson1992airline}
Elizabeth~Louise Williamson.
\newblock \emph{Airline network seat inventory control: Methodologies and
  revenue impacts}.
\newblock PhD thesis, Massachusetts Institute of Technology, 1992.

\bibitem[Ye et~al.(2003)Ye, Yung, and Wang]{ye2003fuzzy}
Cang Ye, Nelson~HC Yung, and Danwei Wang.
\newblock A fuzzy controller with supervised learning assisted reinforcement
  learning algorithm for obstacle avoidance.
\newblock \emph{IEEE Transactions on Systems, Man, and Cybernetics, Part B
  (Cybernetics)}, 33\penalty0 (1):\penalty0 17--27, 2003.

\bibitem[Zaheer et~al.(2017)Zaheer, Kottur, Ravanbakhsh, Poczos, Salakhutdinov,
  and Smola]{zaheer2017deep}
Manzil Zaheer, Satwik Kottur, Siamak Ravanbakhsh, Barnabas Poczos, Russ~R
  Salakhutdinov, and Alexander~J Smola.
\newblock Deep sets.
\newblock In I.~Guyon, U.~Von Luxburg, S.~Bengio, H.~Wallach, R.~Fergus,
  S.~Vishwanathan, and R.~Garnett, editors, \emph{Advances in Neural
  Information Processing Systems}, volume~30. Curran Associates, Inc., 2017.

\end{thebibliography}

\end{document}